\documentclass[11pt,a4paper,twoside,final]{scrartcl}
\usepackage{a4wide}
\usepackage{amsfonts}
\usepackage{amsmath}
\usepackage{amssymb}
\usepackage{amsthm}
\usepackage{mathtools}
\usepackage[mathscr]{euscript}
\usepackage[numbers]{natbib}
\usepackage{booktabs}
\usepackage[Algorithm]{algorithm}
\usepackage{algorithmic}
\usepackage{pgfplots}
\usepackage{pgfplotstable}
\pgfplotsset{compat=newest}
\usepackage{graphicx}
\usepackage{soul}
\usepackage{color}
\usepackage{colortbl}
\usepackage{multirow}
\usepackage{boxedminipage}
\usepackage{enumerate}
\usepackage{caption}
\usepackage{subcaption}
\usepackage{rotating}
\usepackage{hyperref}

\newcommand{\A}{\ensuremath{\mathcal{A}}}
\newcommand{\C}{\ensuremath{\mathbb{C}}}
\newcommand{\D}{\ensuremath{\mathcal{D}}}
\newcommand{\E}{\ensuremath{\mathbb{E}}}
\newcommand{\I}{\ensuremath{\mathrm{I}}}

\newcommand{\N}{\ensuremath{\mathbb{N}}}

\newcommand{\OO}[1]{\mathcal{O}\left(#1\right)}
\renewcommand{\P}{\ensuremath{\mathbb{P}}}

\newcommand{\R}{\ensuremath{\mathbb{R}}}

\newcommand{\T}{\ensuremath{\mathbb{T}}}

\newcommand{\Z}{\ensuremath{\mathbb{Z}}}

\newcommand{\boldf}{{\ensuremath{\boldsymbol{f}}}}
\newcommand{\boldh}{{\ensuremath{\boldsymbol{h}}}}

\newcommand{\boldk}{{\ensuremath{\boldsymbol{k}}}}

\newcommand{\boldu}{{\ensuremath{\boldsymbol{u}}}}

\newcommand{\boldX}{{\ensuremath{\boldsymbol{X}}}}
\newcommand{\boldx}{{\ensuremath{\boldsymbol{x}}}}

\newcommand{\boldgamma}{{\ensuremath{\boldsymbol{\gamma}}}}

\newcommand{\boldxi}{{\ensuremath{\boldsymbol{\xi}}}}

\newcommand{\boldone}{{\ensuremath{\boldsymbol{1}}}}

\newcommand{\comp}{\ensuremath{\mathrm{c}}}
\newcommand{\e}{\textnormal{e}}
\newcommand{\ii}{\textnormal{i}}
\renewcommand{\u}{\ensuremath{\mathfrak{u}}}
\renewcommand{\v}{\ensuremath{\mathfrak{v}}}
\newcommand{\fhat}{{\ensuremath{\hat{f}}}}
\newcommand{\ghat}{{\ensuremath{\hat{g}}}}

\newcommand{\abs}[1]{\left\vert #1\right\vert}
\newcommand{\norm}[1]{\left\Vert #1\right\Vert}
\newcommand{\scalar}[2]{\left\langle #1, #2\right\rangle}

\newcommand{\floor}[1]{\left\lfloor #1\right\rfloor}
\newcommand{\bigtimes}{\mathop{\text{\Large{$\times$}}}}
\renewcommand{\ln}{\mathrm{ln\,}}

\newcommand{\diff}{\ensuremath{\mathrm{d}}}

\DeclareMathOperator*{\esssup}{ess\,sup}

\DeclareMathOperator*{\sinc}{sinc}

\definecolor{light-gray}{gray}{0.75} 
\sethlcolor{light-gray}
\definecolor{color1}{HTML}{E69F00}
\definecolor{color2}{HTML}{56B4E9}
\definecolor{color3}{HTML}{009E73}
\definecolor{color4}{HTML}{F0E442}
\definecolor{color5}{HTML}{0072B2}
\definecolor{color6}{HTML}{D55E00}
\definecolor{color7}{HTML}{CC79A7}

\makeatother

\newtheorem{theorem}{Theorem}[section]
\newtheorem{lemma}[theorem]{Lemma}
\newtheorem{remark}[theorem]{Remark}
\newtheorem{generalisation}[theorem]{Generalisation}
\newtheorem{definition}[theorem]{Definition}
\newtheorem{example}[theorem]{Example}
\newtheorem{corollary}[theorem]{Corollary}
\newtheorem{proposition}[theorem]{Proposition}

\newenvironment{Theorem}{\goodbreak \begin{theorem}\normalfont \slshape}{\end{theorem}}
\newenvironment{Lemma}{\goodbreak \begin{lemma}\normalfont \slshape}{\end{lemma}}
\newenvironment{Remark}{\goodbreak \begin{remark}\normalfont \rmfamily}{\bend\end{remark}}

\newenvironment{Corollary}{\goodbreak \begin{corollary}\normalfont \slshape}{\end{corollary}}

\makeatletter
\def\imod#1{\allowbreak\mkern10mu({\operator@font mod}\,\,#1)}
\makeatother

\numberwithin{equation}{section}
\numberwithin{table}{section}
\numberwithin{figure}{section}

\newcommand{\bend}{\hspace*{0ex} \hfill \hbox{\vrule height
    1.5ex\vbox{\hrule width 1.4ex \vskip 1.4ex\hrule  width 1.4ex}\vrule
    height 1.5ex}}

\long\def\symbolfootnote[#1]#2{\begingroup%
\def\thefootnote{\fnsymbol{footnote}}\footnote[#1]{#2}\endgroup}
\renewcommand{\thefootnote}{\fnsymbol{footnote}}
\title{Nonlinear approximation in bounded orthonormal product bases}

\date{\today}
\author{Lutz K\"ammerer\footnotemark[1], Daniel Potts\footnotemark[2], Fabian Taubert\footnotemark[3]}

\hypersetup{pdfauthor=...,
          pdftitle={Nonlinear approximation in bounded orthonormal product bases},
          pdfsubject={},
          pdfcreator = {pdflatex},
          plainpages=false,
          pdfstartview=FitH,       
          pdfview=FitH,            
          pdfpagemode=UseOutlines, 
          bookmarksnumbered=true, 
          bookmarksopen=false,     
          bookmarksopenlevel=0,   
          colorlinks=true,       
          linkcolor=black,
          citecolor=black,
          urlcolor=black}

\newif\ifshowextendedpaperversion
\showextendedpaperversionfalse

\begin{document}

\maketitle

\begin{abstract}
\small
We present a dimension-incremental algorithm for the nonlinear approximation of high-dimensional functions in an arbitrary bounded orthonormal product basis. Our goal is to detect a suitable truncation of the basis expansion of the function, where the corresponding basis support is assumed to be unknown. Our method is based on point evaluations of the considered function and adaptively builds an index set of a suitable basis support such that the approximately largest basis coefficients are still included. For this purpose, the algorithm only needs a suitable search space that contains the desired index set. Throughout the work, there are various minor modifications of the algorithm discussed as well, which may yield additional benefits in several situations. For the first time, we provide a proof of a detection guarantee for such an index set in the function approximation case under certain assumptions on the sub-methods used within our algorithm, which can be used as a foundation for similar statements in various other situations as well. Some numerical examples in different settings underline the effectiveness and accuracy of our method.

\medskip

{\small%
\noindent {\textit{2020 AMS Mathematics Subject Classification}} : 
%35C09, % Trigonometric solutions to PDEs
%35R60, % PDEs with randomness, stochastic partial differential equations
%60B15, % Probability measures on groups or semigroups, Fourier transforms, factorization
41A50, % Best approximation, Chebyshev systems
42B05, % Fourier series and coefficients in several variables
%42B37, % Harmonic analysis and PDEs
%60-08, % Computational methods for problems pertaining to probability theory
%65C20, % Probabilistic models, generic numerical methods in probability and statistics
%65C30, % Numerical solutions to stochastic differential and integral equations
65D15, % Algorithms for approximation of functions
65D30, % Numerical integration
65D32, % Numerical quadrature and cubature formulas
65D40, % High-dimensional functions; sparse grids
65T40, % Numerical methods for trigonometric approximation and interpolation
%65T50 % Discrete and fast Fourier transforms
65Y20, % Complexity and performance of numerical algorithms
%68Q25, % Analysis of algorithms and problem complexity 
%68W20, % Randomize algorithms
%68W40, % Analysis of algorithms
%94A20  % Sampling theory
}
\end{abstract}

\footnotetext[1]{
  Chemnitz University of Technology, Faculty of Mathematics, 09107 Chemnitz, Germany\\
  kaemmerer@mathematik.tu-chemnitz.de
}
\footnotetext[2]{
  Chemnitz University of Technology, Faculty of Mathematics, 09107 Chemnitz, Germany\\
  potts@mathematik.tu-chemnitz.de
}
\footnotetext[3]{
  Chemnitz University of Technology, Faculty of Mathematics, 09107 Chemnitz, Germany\\
  fabian.taubert@mathematik.tu-chemnitz.de
}
\medskip

\ifshowextendedpaperversion
%\newpage
\tableofcontents
\newpage
\fi

%%%%%%%%%%%%%%%%%%%%%%%%%%%%%%%%%%%%%%%%%%%%%%%%%%%%%%%%%%%%%%%%%%%%%%%%%%%%%%

\section{Introduction}
\label{sec:intro}
In recent years, so-called sparse algorithms that are designed to recover sparse signals have gained significant attention. Various methods and algorithms got developed since then and evolved the field of compressed sensing tremendously, see \cite{FoRa13} for a lot of examples and references. Especially the so-called \textit{sparse Fast Fourier Transform (sFFT)} algorithms (see, e.g., \cite{IwGiSt07,Iw10,HaInKaPr12a,HaInKaPr12,LaWaCh13,SeIw13,Akavia14,ChLaWa15,PlWa16,PlWa16a,Bittens17} or \cite{GiInIwSchm14} for a short overview and introduction) provide efficient ways to reconstruct univariate sparse trigonometric polynomials in different settings. Of course, there are also many other one-dimensional bases besides the Fourier basis, where the sparse recovery problem is also of interest. Hence, similar algorithms began to arise for bases such as, e.g., the Legendre polynomial basis \cite{RaWa12,PePlRo12,PoTa16,HuIwKi17}, as well as for more general settings, cf.\ \cite{GiGuReRuWo20}. Often, the methods of the aforementioned are also applicable in the approximation setting, so if the target function is not sparse itself but assumed to be well approximated by some sparse quantity.

At the same time the high-dimensional generalization of these problems became another topic of research, in particular the question for methods that circumvent the curse of dimensionality (as introduced in \cite{Be61}) to some extent. While again the Fourier case is already studied very well, e.g., \cite{Iw13,InKa14,PoVo14,KaPoVo17,PoVo17,ChChWa19b,ChChWa20,KaKrVo20}, there is very little knowledge about efficient algorithms for other high-dimensional bases.

Sparse polynomial approximation algorithms based on least squares, sparse recovery or compressed sensing have been shown to be quite effective for approximating high-dimensional functions, even in non-Fourier settings, with a relatively small number of samples used. A broad overview of this topic is available in \cite{CoDe15, AdBruWeb17, AdBruWeb22} and the numerous references therein, providing detailed analysis of the strengths and limitations of these methods. One of the main challenges of sparse polynomial approximation is the computational complexity of the matrix-vector multiplications involved in these algorithms. The size of the matrices used therein can grow exponentially with the number of input variables, making the computation time of these algorithms a bottleneck in many applications. For certain problem settings the structure of the matrices or the particular function spaces considered allow for a speedup in computation time since faster matrix-vector multiplication algorithms are available. However, for more general problem settings the computational complexity is still an issue. Recently, more efficient sublinear-time algorithms for bounded orthonormal product bases have been developed in \cite{ChIwKr20,ChIwVo21} and have shown promising results.

Another popular approach in the high-dimensional stochastic setting are sparse polynomial chaos expansions, cf.\ \cite{LueMaSu21} and the references therein. There, random variables are approximated using a subset of the corresponding polynomial orthonormal basis by sparse regression. After each iteration, the candidate subset or the sampling locations are modified until the sparse solution is satisfactory. Note that the concept of sparsity is only used as a tool to find robust solutions in this case and is not the main goal of sparse polynomial chaos expansions. There also exist basis-adaptive sparse PCE approaches as for example described and compared in \cite{LueMaSu22} using and combining various approaches to iteratively build a suitable candidate basis. However, the particular methods have to be chosen carefully since the relative error strongly varies for the different methods. A final model selection after computing several sparse PCE solutions is heavily recommended therein.

Our aim here is the nonlinear approximation of a function $f(\boldx) \coloneqq \sum_{\boldk \in \N^d} c_{\boldk}\Phi_{\boldk}(\boldx)$ using samples by a truncated sum
\begin{align*}
S_\I f(\boldx) \coloneqq \sum_{\boldk \in \I} c_{\boldk}\,\Phi_{\boldk}(\boldx),
\end{align*}
with a carefully selected, finite index set $\I$, which is a-priori unknown. Additionally, we also approximate the basis coefficients $c_\boldk, \boldk \in \I,$ to derive the approximation
\begin{align*}
S_\I f(\boldx) \approx S_\I^\A f(\boldx) \coloneqq \sum_{\boldk \in \I} \fhat_{\boldk}\,\Phi_{\boldk}(\boldx),
\end{align*}
with $\fhat_{\boldk} \in \C, \boldk \in \I$. Throughout the whole paper, sampling is meant w.r.t.\ a black box algorithm that provides the function values $f(\boldx)$ for any sampling node $\boldx$ our algorithm $\A$ requires. Such a black box case can be achieved for example when solving parametric PDEs, see for example \cite{KeWeRi19, KaPoTa22}, where each sample can be computed as the solution of the PDE w.r.t.\ the spatial domain for a fixed random variable. This concept also enables the algorithm to work highly adaptive, since the samples can be suitably chosen in each step, which is not the case when working with given samples. 

We stress on the two parts of this sparse approximation problem, which can be identified by the common error estimate
\begin{align*}
\abs{S_\I^\A f(\boldx)-f(\boldx)} &\leq \abs{S_\I^\A f(\boldx)-S_\I f(\boldx)} + \abs{S_\I f(\boldx) - f(\boldx)}\\
&\leq  \abs{\sum_{\boldk \in \I} (\fhat_{\boldk} - c_{\boldk})\Phi_{\boldk}(\boldx)} + \abs{\sum_{\boldk \in \I^\comp} c_{\boldk}\Phi_{\boldk}(\boldx)}\\
&\leq \sum_{\boldk \in \I} \abs{\fhat_{\boldk} - c_{\boldk}} + \sum_{\boldk \in \I^\comp} \abs{c_{\boldk}},
\end{align*}
where we denote with $\I^\comp \coloneqq \N^d \setminus \I$ the complement of $\I$. We need to compute good approximations $\fhat_{\boldk}$ for the coefficients $c_\boldk$ for all $\boldk\in\I$ to reduce the coefficient approximation error $\sum_{\boldk \in \I} \abs{\fhat_{\boldk} - c_{\boldk}}$ and detect a suitable sparse index set $\I$ such that it contains as many indices $\boldk$, corresponding to the largest coefficients $c_\boldk$ of the function $f$, as possible and therefore minimizing the truncation error $\sum_{\boldk \in \I^\comp} \abs{c_{\boldk}}$. While the coefficient approximation problem for a given index set $\I$ is well-known for many bases, the detection of a good index set $\I$ is rather complicated and will therefore be our primary aim. 

In this paper, we present a dimension-incremental approach for the nonlinear approximation of high-dimensional functions by sparse basis representations applicable to arbitrary bounded orthonormal product basis. The basis indices for these representations are detected adaptively by computing so-called projected coefficients, which are indicating the importance of the corresponding index projections. Our algorithm utilizes suitable methods in the corresponding function spaces to determine sampling nodes and approximate those projected coefficients using the corresponding samples, e.g., by a cubature or least squares method. Therefore, our algorithm benefits tremendously, if those methods are efficient in the sense of sampling or computational complexity.

The paper is organized as follows: In the remaining part of Section \ref{sec:intro}, we briefly introduce the function space setting and explain the concept of projected coefficients. In Section \ref{sec:alg}, we derive our dimension-incremental method and briefly discuss its complexity, the a-priori choice of the search space and alternative increment strategies. Section \ref{sec:theory} contains the derivation of the theoretical main result and Section \ref{sec:num} shows the application of our algorithm to periodic and non-periodic function approximations. Finally, we briefly summarize the results of this work in Section \ref{sec:summary}.

\subsection{Bounded orthonormal product bases}\label{subsec:setup}

We consider $d \in \N$ measure spaces $(\D_j,\A_j,\mu_j), j = 1,\ldots,d,$ with the probability measures $\mu_j$, the $\sigma$-algebras $\A_j$ and the Borel sets $\D_j \subset \R$ for all $j=1,\ldots,d$. As usual, we denote the sets of all functions $f : \D_j \rightarrow \C$ that are square-integrable with respect to $\mu_j$ by $L_2(\D_j,\mu_j)$ and for $\D = \bigtimes_{j=1}^d \D_j \subset \R^d$ the set of all functions $f : \D \rightarrow \C$ that are square-integrable with respect to the product measure $\mu = \bigtimes_{j=1}^d \mu_j$ by $L_2(\D,\mu)$. We assume that the measure spaces $(\D_j,\A_j,\mu_j), j = 1,\ldots,d,$ are such that the $L_2(\D_j,\mu_j), j = 1,\ldots,d,$ are separable Hilbert spaces. Hence, there exists a countable orthonormal basis $\lbrace \phi_{j,k_j}:\D_j\rightarrow\C \,\vert\, k_j \in \N \rbrace$ for each $j = 1,\ldots,d$. Further, the space $L_2(\D,\mu)$ is then also a separable Hilbert space spanned by the orthonormal product basis $\lbrace \Phi_{\boldk}:\D\rightarrow\C \,\vert\, \boldk \in \N^d \rbrace$ with
\begin{align*}
\Phi_{\boldk}(\boldx) \coloneqq \prod_{j=1}^d \phi_{j,k_j}(x_j), && \forall \boldk = (k_j)_{j=1}^d \in \N^d,\, \forall\boldx=(x_j)_{j=1}^d \in \D.
\end{align*}
Finally, we assume that there exist the finite constants
\begin{align*}
B_j \coloneqq \max_{k_j \in \N} \norm{\phi_{j,k_j}}_\infty \coloneqq \max_{k_j \in \N} \esssup_{x_j \in \D_j} \abs{\phi_{j,k_j}(x_j)} \in [1,\infty)
\end{align*}
for each $j=1,\ldots,d$, i.e., the orthonormal basis $\lbrace \phi_{j,k_j}:\D_j\rightarrow\C \,\vert\, k_j \in \N \rbrace$ is bounded for each $j$. Then, the orthonormal product basis $\lbrace \Phi_{\boldk}:\D\rightarrow\C \,\vert\, \boldk \in \N^d \rbrace$ is also bounded by
\begin{align*}
B \coloneqq \prod_{j=1}^d B_j = \max_{\boldk \in \N^d} \norm{\Phi_{\boldk}}_\infty = \max_{\boldk \in \N^d} \esssup_{\boldx \in \D} \abs{\Phi_{\boldk}(\boldx)} \in [1,\infty)
\end{align*}
and is therefore called \textit{Bounded Orthonormal Product Basis (BOPB)} throughout this paper.

Let $f \in L_2(\D,\mu)$ be smooth enough such that there exist coefficients $\lbrace c_{\boldk}\rbrace_{\boldk \in \N^d}, \sum_{\boldk \in \N^d} \abs{c_{\boldk}} < \infty,$ and the series expansion
\begin{align*}
f(\boldx) \coloneqq \sum_{\boldk \in \N^d} c_{\boldk}\Phi_{\boldk}(\boldx)
\end{align*}
holds pointwise for all $\boldx \in \D$. This smoothness requirement enables the concept of approximation of $f$ using point samples, but is different for each basis $\lbrace \Phi_{\boldk}\rbrace_{\boldk \in \N^d}$.

An example of such a BOPB is the Fourier system on the periodic, $d$-variate torus $\D = \T^d \simeq [0,1)^d$ with the common Lebesgue measure, where $\Phi_\boldk (\boldx) = \e^{2\pi\ii\boldk\cdot\boldx}, \boldk \in \Z^d,$ with constant $B=1$. For $f \in L_2(\T^d)$, the coefficients $c_\boldk$ are then the well known Fourier coefficients $\int_{\T^d} f(\boldx) \e^{-2\pi\ii\boldk\cdot\boldx} \diff \boldx$. The approximation of a Fourier partial sum $S_\I f(\boldx)$ for a given frequency set $\I$ can be realized efficiently by specific Fast Fourier Transform (FFT) methods, while the more challenging task of identifying a suitable frequency set $\I$ is considered as ``sparse FFT'' in several works, see \cite[Tbl. 1.1]{KaKrVo20} for an overview.

Another example is the Chebyshev system on $\D = [-1,1]^d$ with $\Phi_\boldk(\boldx) = T_\boldk(\boldx) \coloneqq \prod_{j=1}^d \cos(k_j \arccos(x_j)), \boldk \in \N^d,$ a tensor product of Chebyshev polynomials of first kind. The corresponding space is $L_2([-1,1]^d,\mu_{\mathrm{Cheb}})$, where $\mu_{\mathrm{Cheb}}(\diff\boldx) \coloneqq \prod_{j=1}^d (\pi\sqrt{1-x_j^2})^{-1} \diff\boldx$ is the Chebyshev measure. The BOPB constant is again $B=1$ in this setting. 

We encourage the reader to keep such an example in mind. Especially for the Fourier system, the sparse FFT approaches presented in \cite{PoVo14,KaPoVo17,KaKrVo20} may be seen as special cases of the algorithm we are about to present. Note that our setting above is not restricted to bases with similar structures in each dimension $j=1,\ldots,d$. For instance, one could also think of systems with a Fourier-type basis for only some $j \in \{1,\ldots,d\}$ and a Chebyshev-type basis for the remaining dimensions, see \cite[Sec.\ 5.1.1]{ChIwVo21} as an example.

\begin{Remark}\label{rem:smooth}
As mentioned above, the smoothness of $f$ is important to ensure the well-definedness of point evaluations of $f$. Obviously, a possible restriction to ensure this property is to assume the continuity of $f$. Additionally, the smoothness condition is also fulfilled for most function spaces with higher regularity, e.g., when using weighted Wiener spaces as considered in \cite{PoSchmi19,JaUlVo22} for the Fourier setting. However, depending on the BOPB, other or weaker assumptions may be possible.

Also, one could assume $f$ to be from a reproducing kernel Hilbert space instead, cf.\ for example \cite{BeTh04}. The well-definedness of point evaluations is one of the defining properties of such spaces.
\end{Remark}

\subsection{Projected coefficients and cubature rules}\label{subsec:pc}

To simplify notations in the upcoming sections, we introduce the following notations for $\u \subset \{1,\ldots,d\}$ and its complement $\u^\comp \coloneqq \{1,\ldots,d\} \setminus \u$:
\begin{itemize}
\item $\D_{\u} \coloneqq \bigtimes_{j \in \u} \D_j \subset \R^{\abs{\u}}$, $\mu_{\u} \coloneqq \bigtimes_{j \in \u} \mu_j$, $B_\u \coloneqq \prod_{j\in\u} B_j$,
\item $\Phi_{\u,\boldk}(\boldxi) \coloneqq \prod_{j \in \u} \phi_{j,k_j}(\xi_j)$ for all $\boldk=(k_j)_{j \in \u}\in\N^{\abs{\u}}$ and $\boldxi=(\xi_j)_{j\in\u}\in\D_{\u}$,
\item $\boldh_{\u} \coloneqq (h_j)_{j\in\u} \in \N^{\abs{\u}}$ for $\boldh \in \N^d$,
\item $(\boldk,\boldh)_{\u} \coloneqq (l_j)_{j=1}^d$ with $l_j = \begin{cases} k_j, & j \in \u \\ h_j, & j \not\in\u  \end{cases}$ for $\boldk = (k_j)_{j\in\u} \in \N^{\abs{\u}},\;\boldh = (h_j)_{j\in\u^\comp} \in \N^{d-\abs{\u}}$,
\item $f(\boldxi,\tilde{\boldx})_{\u} \coloneqq f((y_j)_{j=1}^d)$ with $y_j = \begin{cases} \xi_j, & j \in \u \\ \tilde{x}_j, & j \not\in\u  \end{cases}$ for $\boldxi = (\xi_j)_{j\in\u} \in \D_{\u},\;\tilde{\boldx} = (\tilde{x}_j)_{j\in\u^\comp} \in \D_{\u^\comp}$.
\end{itemize}
To ensure that all of these quantities using $\u$ (or $\u^\comp$) are well-defined, we assume them to be ordered naturally. Note that the notations coincide with their one- and $d$-dimensional counterparts for $\abs{\u}=1$ or $\abs{\u}=d$, respectively, if they exist.

Our algorithm, which we are about to present in Section \ref{sec:alg}, detects a suitable index set $\I$ by computing approximations of so-called projected coefficients $c_{\u,\boldk}$ for $\u = \{1,\ldots,t\}$ and $\u = \{t\}$ for each $t = 1,\ldots,d$ and several indices $\boldk \in \N^{\abs{\u}}$ using samples of the function $f$. In particular, we consider the projected coefficients
\begin{align*}
c_{\u,\boldk}(\tilde{\boldx}) \coloneqq \scalar{f(\cdot,\tilde{\boldx})_{\u}}{\Phi_{\u,\boldk}(\cdot)}_{\D_{\u}} = \int_{\D_{\u}} f(\boldxi,\tilde{\boldx})_{\u} \,\overline{\Phi_{\u,\boldk}(\boldxi)}\,\diff\mu_{\u}(\boldxi)
\end{align*}
as a function w.r.t.\ the $(d-\abs{\u})$-dimensional anchor $\tilde{\boldx}\in\D_{\u^\comp}$. The name is based on the fact that those $c_{\u,\boldk}(\tilde{\boldx})$ can be interpreted as the coefficients of the basis expansion in the space $L_2(\D_\u,\mu_\u)$ of the projections $f(\cdot,\tilde{\boldx})$, which play an important role in the anchored version of the multivariate decomposition method (MDM), cf.\ for example \cite{KuSloWaWo10}.

However, using $\sum_{\boldh \in \N^d} \abs{c_{\boldh}} < \infty$ and Fubini's Theorem, we proceed
\begin{align}
c_{\u,\boldk}(\tilde{\boldx}) &= \int_{\D_{\u}} f(\boldxi,\tilde{\boldx})_{\u} \,\overline{\Phi_{\u,\boldk}(\boldxi)}\,\diff\mu_{\u}(\boldxi) \label{eq:proj_coef} \\
&= \int_{\D_{\u}} \left( \sum_{\boldh\in\N^d} c_{\boldh} \Phi_{\boldh}(\boldxi,\tilde{\boldx})_{\u} \right) \,\overline{\Phi_{\u,\boldk}(\boldxi)}\,\diff\mu_{\u}(\boldxi) \nonumber  \\
&= \int_{\D_{\u}} \sum_{\boldh\in\N^d} c_{\boldh} \Phi_{\u,\boldh_{\u}}(\boldxi)\Phi_{\u^\comp,\boldh_{\u^\comp}}(\tilde{\boldx}) \,\overline{\Phi_{\u,\boldk}(\boldxi)}\,\diff\mu_{\u}(\boldxi) \nonumber \\
&= \sum_{\boldh\in\N^d} c_{\boldh} \Phi_{\u^\comp,\boldh_{\u^\comp}}(\tilde{\boldx}) \int_{\D_{\u}} \Phi_{\u,\boldh_{\u}}(\boldxi) \,\overline{\Phi_{\u,\boldk}(\boldxi)}\,\diff\mu_{\u}(\boldxi) \nonumber \\
&= \sum_{\boldh\in\N^d} c_{\boldh} \Phi_{\u^\comp,\boldh_{\u^\comp}}(\tilde{\boldx}) \scalar{\Phi_{\u,\boldh_{\u}}}{\Phi_{\u,\boldk}}_{\D_{\u}} \nonumber \\
&= \sum_{\boldh=(\boldk,\boldh_{\u^\comp})_{\u}\in\N^d} c_{\boldh} \Phi_{\u^\comp,\boldh_{\u^\comp}}(\tilde{\boldx}). \label{eq:proj_coef_simpl}
\end{align}
Hence, the size of the projected coefficients $c_{\u,\boldk}$ can be considered as an indicator for the importance of the set of indices $\boldh=(\boldk,\boldh_{\u^\comp})_{\u}$ with fixed $\boldk$ in the components $\u$.

In order to utilize this fact, we need a suitable way to approximate such projected coefficients. Here, one can apply various approaches, which we will call reconstruction method throughout this paper. For our theoretical results, we will restrict ourselves to a special kind of cubature approaches in the following sections. Note that most of the theoretical results in Section \ref{sec:theory} can be proven similarly for other reconstruction methods, cf.\ Remark \ref{rem:lsqr}.

We require for fixed $\u$ a suitable cubature rule $Q$ with weights $w_j \in \C, j=1,\ldots,M$ and cubature nodes $\boldxi_j \in \D_{\u},j=1,\ldots,M$, which is exact for some finite index set $K \subset \N^{\abs{\u}}$ for the inner products $\scalar{\Phi_{\u,\boldk_1}}{\Phi_{\u,\boldk_2}}_{\D_{\u}}$ for all $\boldk_1,\boldk_2 \in K$, i.e.,
\begin{align}\label{eq:quadrature_exact}
\sum_{j=1}^M w_j \Phi_{\u,\boldk_1}(\boldxi_j) \,\overline{\Phi_{\u,\boldk_2}(\boldxi_j)} = \delta_{\boldk_1,\boldk_2} = \begin{cases} 1 & \text{if $\boldk_1 = \boldk_2$} \\
0 & \text{else}
\end{cases} && \forall \boldk_1,\boldk_2 \in K
\end{align}
holds. Additionally, we denote
\begin{align}\label{eq:weight_sum}
C_Q \coloneqq \sum_{j=1}^M \abs{w_j}
\end{align} 
for each such cubature rule.

We now define the approximated projected coefficients with anchor $\tilde{\boldx}$ as cubature of the integral \eqref{eq:proj_coef} w.r.t.\ the cubature rule $Q$, i.e.,
\begin{align}\label{eq:proj_coef_approx}
\fhat_{\u,\boldk}^Q (\tilde{\boldx}) \coloneqq \sum_{j=1}^M w_j f(\boldxi_j,\tilde{\boldx})_\u \overline{\Phi_{\u,\boldk}(\boldxi_j)}.
\end{align}
Note that the approximation of the projected coefficients $c_{\u,\boldk}(\tilde{\boldx})$ with anchor $\tilde{\boldx}$ may also be realized in different ways, cf.\ Remark \ref{rem:lsqr}.

With similar arguments as above, we get
\begin{align*}
\fhat_{\u,\boldk}^Q (\tilde{\boldx}) &= \sum_{j=1}^M w_j f(\boldxi_j,\tilde{\boldx})_{\u}\overline{\Phi_{\u,\boldk}(\boldxi_j)} \nonumber \\
&= \sum_{j=1}^M w_j \left( \sum_{\boldh\in\N^d} c_{\boldh} \Phi_{\boldh}(\boldxi_j,\tilde{\boldx})_\u \right) \,\overline{\Phi_{\u,\boldk}(\boldxi_j)} \nonumber \\
&= \sum_{\boldh\in\N^d} c_{\boldh} \Phi_{\u^\comp,\boldh_{\u^\comp}}(\tilde{\boldx}) \sum_{j=1}^M w_j \Phi_{\u,\boldh_\u}(\boldxi_j) \,\overline{\Phi_{\u,\boldk}(\boldxi_j)}.
\end{align*}
We assume now that \eqref{eq:quadrature_exact} holds for some index set $K \subset \N^{\abs{\u}}$ and consider another index set $J \subset \N^d$ with $J \subset K \times \N^{d-\abs{\u}}$. We split the sum $\sum_{\boldh \in \N^d} = \sum_{\boldh \in J} + \sum_{\boldh \in J^\comp}$, apply \eqref{eq:quadrature_exact} in the first sum $\sum_{\boldh \in J}$ and continue for all $\boldk \in K$ with
\begin{align}\label{eq:proj_approx_early}
\fhat_{\u,\boldk}^Q (\tilde{\boldx}) = &\sum_{\boldh = (\boldk,\boldh_{\u^\comp})_\u\in J} c_{\boldh} \Phi_{\u^\comp,\boldh_{\u^\comp}}(\tilde{\boldx}) +\sum_{\boldh\in J^\comp} c_{\boldh} \Phi_{\u^\comp,\boldh_{\u^\comp}}(\tilde{\boldx}) \sum_{j=1}^M w_j \Phi_{\u,\boldh_\u}(\boldxi_j) \,\overline{\Phi_{\u,\boldk}(\boldxi_j)},
\end{align}
which is the same as \eqref{eq:proj_coef_simpl} up to the projection error term
\begin{align}\label{eq:noise}
\Psi_{\u,\boldk}^{Q,J}(\tilde{\boldx}) \coloneqq \sum_{\boldh\in J^\comp} c_{\boldh} \Phi_{\u^\comp,\boldh_{\u^\comp}}(\tilde{\boldx}) \sum_{j=1}^M w_j \Phi_{\u,\boldh_\u}(\boldxi_j) \,\overline{\Phi_{\u,\boldk}(\boldxi_j)}.
\end{align}
Note that $\Psi_{\u,\boldk}^{Q,J}$ vanishes for sparse functions $f$, i.e., if all the coefficients $c_{\boldh}, \boldh \in J^\comp,$ are zero. Formula \eqref{eq:proj_approx_early} legitimizes the use of $\fhat_{\u,\boldk}^Q$ instead of $c_{\u,\boldk}$ as an indicator for the importance of the $\boldh=(\boldk,\boldh_{\u^\comp})_\u$, if the projection error term $\Psi_{\u,\boldk}^{Q,J}$ is suitably bounded, cf.\ Section \ref{sec:theory}.

\begin{Remark}
Our exactness condition \eqref{eq:quadrature_exact} can be extended to hold for all functions in the span of the respective basis functions $\Phi_{\boldu,\boldk}$ by linearity. It is shown in \cite[Thm.~2.3]{BaKaPoUl22}, that such a condition is equivalent to fulfulling an $L_2$-Marcinkiewicz-Zygmund inequality with equal constants $A=B$.

Unfortunately, this special kind of $L_2$-MZ inequality does only hold for one of our used reconstruction methods in Section \ref{sec:num}, namely the single rank-1 lattice (R1L) approach. Based on this observation we assume, that a generalization of the theoretical part in Section \ref{sec:theory} is possible when assuming the reconstruction method $Q$ to fulfill a relaxed version of the $L_2$-MZ inequality with constants $A \leq B$. Such a condition also holds for the Monte Carlo nodes (MC, cMC) and probably even for the multiple rank-1 lattice (MR1L, cMR1L) approaches from our numerical tests in Section \ref{sec:num}.
\end{Remark}

\begin{Remark}\label{rem:lsqr}
We can write \eqref{eq:proj_coef_approx} as matrix vector equation $\hat{\boldf}_\u^Q = \boldsymbol{\Phi}_\u^\ast \boldf_w$ with $\boldf_w \coloneqq (w_j f(\boldxi_j,\tilde{\boldx})_\u)_{j=1}^M$ and $\boldsymbol{\Phi}_\u \in \C^{M\times \vert K \vert}$ containing the corresponding basis function values $\Phi_{\u,\boldk}(\boldxi_j)$. While we stick to the presented cubature approach $Q$ in the theoretical part of this paper, one can also apply other reconstruction approaches $R$ to compute the approximated projected coefficients $\hat{\boldf}_\u^R$, e.g., using a least squares or compressed sensing approach, cf. \cite[Chap.~3]{FoRa13} for some basic methods. Then, the approximated projected coefficients $\fhat_{\u,\boldk}^{R} (\tilde{\boldx})$ are still a good indicator for the importance of the corresponding indices $\boldh=(\boldk,\boldh_{\u^\comp})_\u$ as long as the corresponding projection error term $\abs{\fhat_{\u,\boldk}^{R} (\tilde{\boldx}) - c_{\u,\boldk}(\tilde{\boldx})}$ is small enough.

Note that the theoretical results studied in Section \ref{sec:theory} should be applied for the new projection error term $\fhat_{\u,\boldk}^{R} (\tilde{\boldx}) - c_{\u,\boldk}(\tilde{\boldx})$ instead of $\Psi_{\u,\boldk}^{Q,J}(\tilde{\boldx})$ in this case and may need some modifications based on the properties of this new projection error term.
\end{Remark}

\iffalse
\begin{Remark}\label{rem:lsqr}
While we stick to the presented cubature approach $Q$ in the theoretical part of this paper, one can also use other reconstruction approaches $R$, e.g., a least squares method, to compute the approximation of the projected coefficients $c_{\u,\boldk}(\tilde{\boldx})$ with anchor $\tilde{\boldx}$, which we also utilize in our numerical experiments in Section \ref{sec:num}. Then, the approximated projected coefficients $\fhat_{\u,\boldk}^{R} (\tilde{\boldx})$ are still a good indicator for the importance of the corresponding indices $\boldh=(\boldk,\boldh_{\u^\comp})_\u$ as long as the corresponding projection error term $\abs{\fhat_{\u,\boldk}^{R} (\tilde{\boldx}) - c_{\u,\boldk}(\tilde{\boldx})}$ is small enough.

Note that the theoretical results studied in Section \ref{sec:theory} should be applied for the new projection error term $\fhat_{\u,\boldk}^{R} (\tilde{\boldx}) - c_{\u,\boldk}(\tilde{\boldx})$ instead of $\Psi_{\u,\boldk}^{Q,J}(\tilde{\boldx})$ in this case and may need some modifications based on the properties of this new projection error term.
\end{Remark}
\fi

\section{The nonlinear approximation algorithm}\label{sec:alg}

In this section, we present our nonlinear approximation algorithm based on the concept of projected coefficients explained in Section \ref{subsec:pc}. In Section \ref{subsec:alternative} we also discuss different increment strategies and their possible advantages and disadvantages.

\subsection{The dimension-incremental method}
\label{subsec:method}

The full method is given in Algorithm \ref{alg:main}. Additionally, Figure \ref{fig:algo} illustrates some of the first steps of the application of Algorithm \ref{alg:main} to some made up function $f$.

As already mentioned, our algorithm proceeds in a dimension-incremental way. Roughly speaking, it constructs the frequencies $\boldk$ of the desired index set $\I$ component-by-component. To explain this concept properly, we denote with $\mathscr{P}(\Omega)$ the power set of a set $\Omega$ and introduce the projection operator 
\begin{align*}
\mathcal{P}_{\u}: \mathscr{P}(\N^d) &\rightarrow \mathscr{P}(\N^{\abs{\u}})\\
\mathcal{P}_{\u}(\Omega) &= \left\lbrace \boldk \in \N^{\abs{\u}} \,\vert\, \exists\, \boldh \in \Omega: \boldk = \boldh_{\u} \right\rbrace.
\end{align*}
Hence, the set $\mathcal{P}_{\u}(\Omega)$ contains all indices $\boldk$ which can be extended to at least one index $\boldh \in \Omega$, i.e., $\boldh = (\boldk,\boldh_{\u^\comp})_\u$ for some $\boldh \in \Omega$.

\subsubsection*{Single component identification}

Algorithm \ref{alg:main} starts by detecting one-dimensional index sets, which we denote as $\I_{\{t\}}$, for all $t=1,\ldots,d$ in step 1. To this end, it constructs a suitable cubature rule $Q$ to compute the approximated projected coefficients $\fhat_{\{t\},k_t}^Q(\tilde{\boldx})$ for all $k_t \in \mathcal{P}_{\{t\}}(\Gamma)$, so all possible values for the $t$-th component of the indices in $\I$ according to our search space $\Gamma$, for some randomly chosen anchor $\tilde{\boldx} \in \D_{\{t\}^\comp}$ via \eqref{eq:proj_coef_approx}. As explained in Section \ref{subsec:pc}, these values are a suitable indicator to decide, whether or not $k_t$ appears as $t$-th component of any index $\boldh \in \I$, i.e., if $k_t \in \mathcal{P}_{\{t\}}(\I)$. An index $k_t$ is kept and therefore added to the set $\I_{\{t\}}$, if the absolute value of its approximated projected  coefficient $\fhat_{\{t\},k_t}^Q(\tilde{\boldx})$ is larger than the so-called detection threshold $\delta_+ \in \R^+$, as it can be seen in Figure \ref{subfig:a}. The right choice of $\delta_+$ and the connection between the detection threshold $\delta_+$ and the true size of the basis coefficients $c_\boldh$ with $\boldh \in \I$ is given in Section \ref{sec:theory}. To avoid the detection of unnecessarily many indices $k_t$, we also use a so-called sparsity parameter $s_\mathrm{local}$ and consider only the $s_\mathrm{local}$-largest approximated projected coefficients $\fhat_{\{t\},k_t}^Q(\tilde{\boldx})$ larger than the detection threshold $\delta_+$. Finally, the random choice of the anchor $\tilde{\boldx}$ may result in some annihilations, so small approximated projected coefficients $\fhat_{\{t\},k_t}^Q(\tilde{\boldx})$ even though the corresponding basis coefficients $c_\boldh$ with $\boldh_{\{t\}} = k_t,$ are large. Therefore, we repeat the choice of $\tilde{\boldx}$, the computation of the approximated projected coefficients $\fhat_{\{t\},k_t}^Q(\tilde{\boldx})$ and the addition of important indices $k_t$ to the index set $\I_{\{t\}}$ now $r \in \N$ times, which is also shown in Figure \ref{subfig:a}. Hence, we call the parameter $r$ the number of detection iterations. Choosing $r$ large enough, cf.\ Section \ref{sec:theory}, ensures that each index $k_t \in \mathcal{P}_{\{t\}}(\I)$ is detected in at least one detection iteration with high probability and therefore $\mathcal{P}_{\{t\}}(\I) \subset \I_{\{t\}}$.

\subsubsection*{Coupled component identification}

In each iteration $t=2,\ldots,d$ of step 2 of Algorithm \ref{alg:main}, we already detected the previous set $\I_{\{1,\ldots,t-1\}}$ and consider $\u = \{1,\ldots,t\}$. As in step 1, the aim is the construction of an index set $\I_\u$ such that $\mathcal{P}_\u(\I) \subset \I_\u$ holds with high probability. We construct our so-called candidate set $K \supset \mathcal{P}_\u(\I)$ from two parts. The first part is the product set $\I_{\{1,\ldots,t-1\}} \times \I_{\{t\}}$. The first set hopefully contains all $\boldk \in \mathcal{P}_{\{1,\ldots,t-1\}}(\I)$ and the second set all $k_t \in \mathcal{P}_{\{t\}}(\I)$. Hence, the combined set $\I_{\{1,\ldots,t-1\}} \times \I_{\{t\}}$ is an obvious choice when we are looking for indices $\boldk \in \mathcal{P}_\u(\I)$. Two such combined sets are shown in Figure \ref{subfig:b} and \ref{subfig:d}. The second part is the projection $\mathcal{P}_\u(\Gamma)$ as in step 1. There is no need to consider any $\boldk \not \in \mathcal{P}_\u(\Gamma)$, since $(\boldk,\boldh)_\u \not \in \Gamma \supset \I$ for any $\boldh \in \N^{d-t}$ anyway in this case. Therefore, the candidate set $K$ is now chosen as the intersection of those two sets, i.e., $K = \left(\I_{\{1,\ldots,t-1\}} \times \I_{\{t\}}\right) \cap \mathcal{P}_\u(\Gamma)$.

Now, we construct a suitable cubature rule $Q$ for the set $K$ and proceed as in the first step of Algorithm \ref{alg:main}: We choose an anchor $\tilde{\boldx} \in \D_{\u^\comp}$ at random, compute the corresponding approximated projected coefficient $\fhat_{\u,\boldk}^Q(\tilde{\boldx})$ and put the indices $\boldk$ of the (up to) $s_\mathrm{local}$-largest coefficients, which are still larger than the detection threshold $\delta_+$, into the set $\I_\u$. For $t=2$, this step is illustrated in Figure \ref{subfig:c}. Finally, we again repeat this procedure $r$ times to ensure the detection of all of the desired indices with high probability. Note that in the final iteration $t=d$ no more than one detection iteration is needed, since the cubature nodes $\boldxi_j$ are already $d$-dimensional, so there is no randomly chosen anchor $\tilde{\boldx}$. Because of this and since the output $\I_\u$ of this final iteration is also the final output $\I$, one might want to use another, smaller sparsity parameter $s$ than in the previous steps. Finally note that the computed approximated projected coefficients in this step are already approximations of the true coefficients $c_\boldh, \boldh \in \I$, so it is not necessary to recompute those quantities in step 3 of Algorithm \ref{alg:main}.

\begin{algorithm}[t]
	\caption{Dimension-incremental Algorithm}\label{alg:main}
  \begin{small}
	\begin{tabular}{p{1.2cm}p{2.3cm}p{10.8cm}}
		Input:  %
		& $\Gamma\subset\N^d$ \hfill & search space \\
		& $f$ & function $f$ as black box (function handle) \\
		& $s,s_\mathrm{local}\in\N$ & sparsity parameter, $s\leq s_\mathrm{local}$ \\  
		& $\delta_+\in\R^+$ & detection threshold \\
		& $r\in\N$ & number of detection iterations%
	\end{tabular}
	\begin{algorithmic}
		\item[(Step 1)] [Single component identification]
		\STATE {\bfseries for} $t:=1,\ldots,d$ {\bfseries do}
		\STATE \hspace{1em} Set $\u = \{t\}$ and $\I_{\u}:=\emptyset$.
		\STATE \hspace{1em} Find $M,$ cubature weights $w_j \in \R,$ and cubature nodes $\boldxi_j \in \D_{\u}$ for $j=1,\ldots,M$ such that \eqref{eq:quadrature_exact} holds for $\mathcal{P}_{\u}(\Gamma)$.
		\STATE \hspace{1em} {\bfseries for} $i:=1,\ldots,r$ {\bfseries do}
		\STATE \hspace{2em} Choose $\tilde{\boldx}\in\D_{\u^\comp}$ at random w.r.t.\ the product measure $\mu_{\u^\comp}$.
		\STATE \hspace{2em} Sample $f$ at the nodes $\left((\boldxi_j,\tilde{\boldx})_\u\right)_{j=1}^M$.
		\STATE \hspace{2em} Compute $\left(\fhat_{\u,k_t}^Q(\tilde{\boldx})\right)_{k_t \in \mathcal{P}_\u(\Gamma)}$ efficiently, cf. \eqref{eq:proj_coef_approx}.
  		\STATE \hspace{2em} Set $\I_\u:=\I_\u\cup\left\lbrace k_t\in\mathcal{P}_\u(\Gamma)\colon \text{(up to) } s_\mathrm{local} \text{-largest values } \abs{\fhat_{\u,k_t}^Q(\tilde{\boldx})} \geq \delta_+\right\rbrace$
		\STATE \hspace{1em} {\bfseries end for} $i$
		\STATE {\bfseries end for} $t$
		\item[(Step 2)] [Coupled component identification]
		\STATE {\bfseries for} $t:=2,\ldots,d$ {\bfseries do}
		\STATE \hspace{1em} If $t<d$, set $\tilde r:=r$ and $\tilde s:=s_\mathrm{local}$, otherwise $\tilde r:=1$ and $\tilde s:=s$. 
		\STATE \hspace{1em} Set $\u = \{1,\ldots,t\}$ and $\I_{\u}:=\emptyset$.
		\STATE \hspace{1em} Construct the index set $K:=(\I_{\{1,\ldots,t-1\}} \times \I_{\{t\}})\cap\mathcal{P}_{\u}(\Gamma)$.
		\STATE \hspace{1em} Find $M,$ cubature weights $w_j \in \R,$ and cubature nodes $\boldxi_j \in \D_{\u}$ for $j=1,\ldots,M$ such that \eqref{eq:quadrature_exact} holds for $K$.
		\STATE \hspace{1em} {\bfseries for} $i:=1,\ldots,\tilde r$ {\bfseries do}
		\STATE \hspace{2em} Choose $\tilde{\boldx}\in\D_{\u^\comp}$ at random w.r.t.\ the product measure $\mu_{\u^\comp}$.
		\STATE \hspace{2em} Sample $f$ at the nodes $\left((\boldxi_j,\tilde{\boldx})_\u\right)_{j=1}^M$.
		\STATE \hspace{2em} Compute $\left( \fhat_{\u,\boldk}^Q(\tilde{\boldx})\right)_{\boldk \in K}$ efficiently, cf. \eqref{eq:proj_coef_approx}.
		\STATE \hspace{2em} Set $\I_{\u}:=\I_{\u}\cup \left\lbrace \boldk \in K\colon \text{(up to) } s_\mathrm{local} \text{-largest values } \abs{\fhat_{\u,\boldk}^Q(\tilde{\boldx})} \geq \delta_+\right\rbrace$.
		\STATE \hspace{1em} {\bfseries end for} $i$
		\STATE {\bfseries end for} $t$
		\item[(Step 3)] %[Computation of basis coefficients]
		%
		%\STATE Generate a sampling set $\mathcal{X}\subset\D$ for $\I^{(1,\ldots,d)}$ such that the corresponding Fourier matrix $\boldsymbol{A}(\mathcal{X},\I^{(1,\ldots,d)})$ is of full column rank and its pseudoinverse can be applied {\bfseries efficiently}.
		%
		%\STATE Compute the basis coefficients $\left(\fhat_{(1,\ldots,d),\boldk}^Q\right)_{\boldk\in \I^{(1,\ldots,d)}}$.
		%
		\STATE Set $\I:=\I_{\u}$ and $\left(\fhat_\boldk\right)_{\boldk\in \I_{\u}} \coloneqq \left(\fhat_{\u,\boldk}^Q\right)_{\boldk\in \I_{\u}}$ with $\u = \{1,\ldots,d\}$.
	\end{algorithmic}

	\begin{tabular}{p{1.3cm}p{2.3cm}p{10.7cm}}
		Output: & $\I\subset\Gamma\subset\N^d$ & detected index set\\
		& $\left(\fhat_\boldk\right)_{\boldk\in \I} \in\C^{\abs{\I}}$ & approximated coefficients with $\abs{\fhat_{\boldk}} \geq \delta_+$ \\
	\end{tabular}
  \end{small}
\end{algorithm}

\begin{figure}
\begin{subfigure}[c]{0.95\textwidth}
	\begin{center}
		\includegraphics{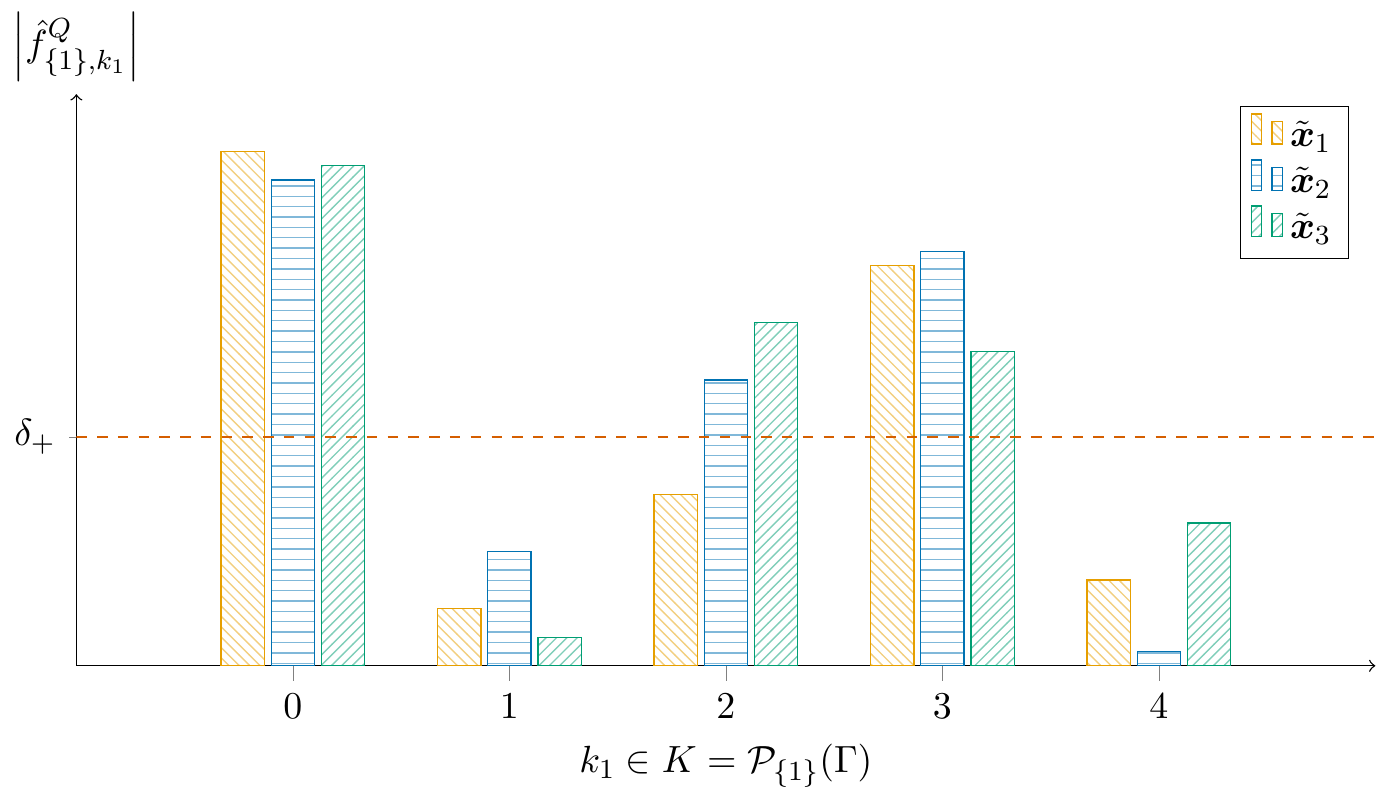}		
	\end{center}
	\caption{One-dimensional detection step: The bars show the absolute values of the approximated projected coefficients $\fhat_{\lbrace 1\rbrace, k_1}^Q$ for $r=3$ random anchors $\tilde{\boldx}$ for all indices $k_1$ in the first candidate set $K=\mathcal{P}_{\{1\}}(\Gamma)$. For $k_1=1$ and $k_1=4$, they are smaller than the detection threshold $\delta_+$ in all $r=3$ detection iterations. Hence, only the indices $0$, $2$ and $3$ are part of $\I_{\lbrace 1 \rbrace}$.}
	\label{subfig:a}
\end{subfigure}
~
\begin{subfigure}[c]{0.95\textwidth}
	\begin{center}
		\includegraphics{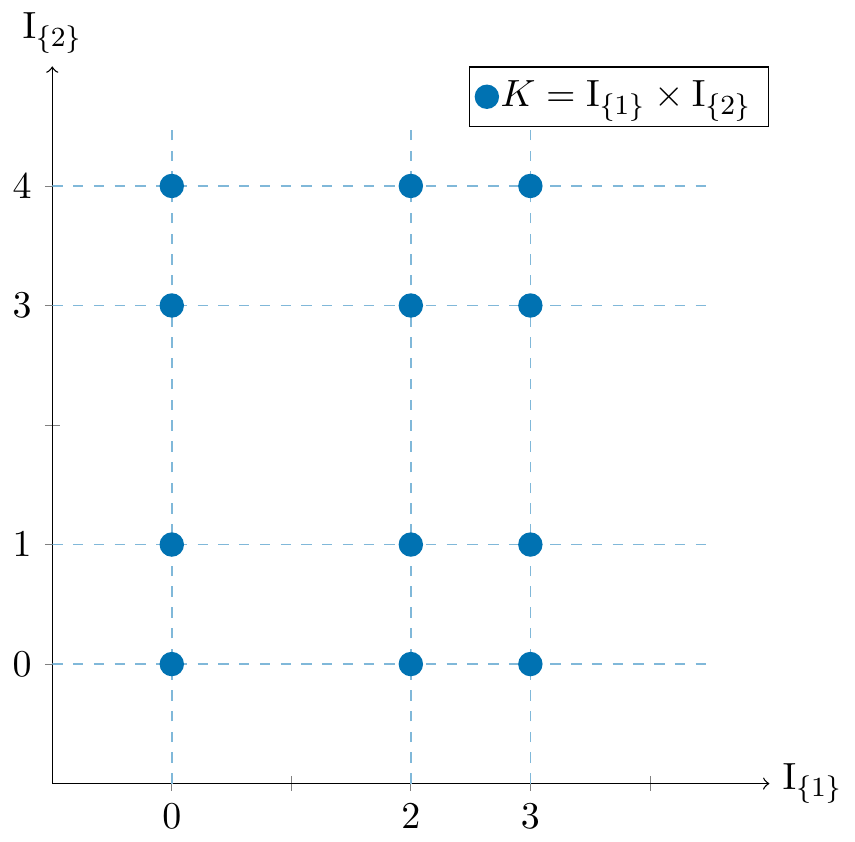}		
	\end{center}
	\caption{Two-dimensional candidate set: The new candidate set $K$ is now just the cartesian product of the one-dimensional detected sets $\I_{\lbrace 1 \rbrace} = \{0,2,3\}$ and $\I_{\lbrace 2 \rbrace}=\{0,1,3,4\}$. The intersection with $\mathcal{P}_{\{1,2\}}(\Gamma)$ is not necessary here, since we assumed $\Gamma$ to be the full grid from $0$ to $4$ in $d$ dimensions.}
	\label{subfig:b}
\end{subfigure}
\end{figure}
\begin{figure}
\ContinuedFloat
\begin{subfigure}[c]{0.95\textwidth}
	\begin{center}
		\includegraphics{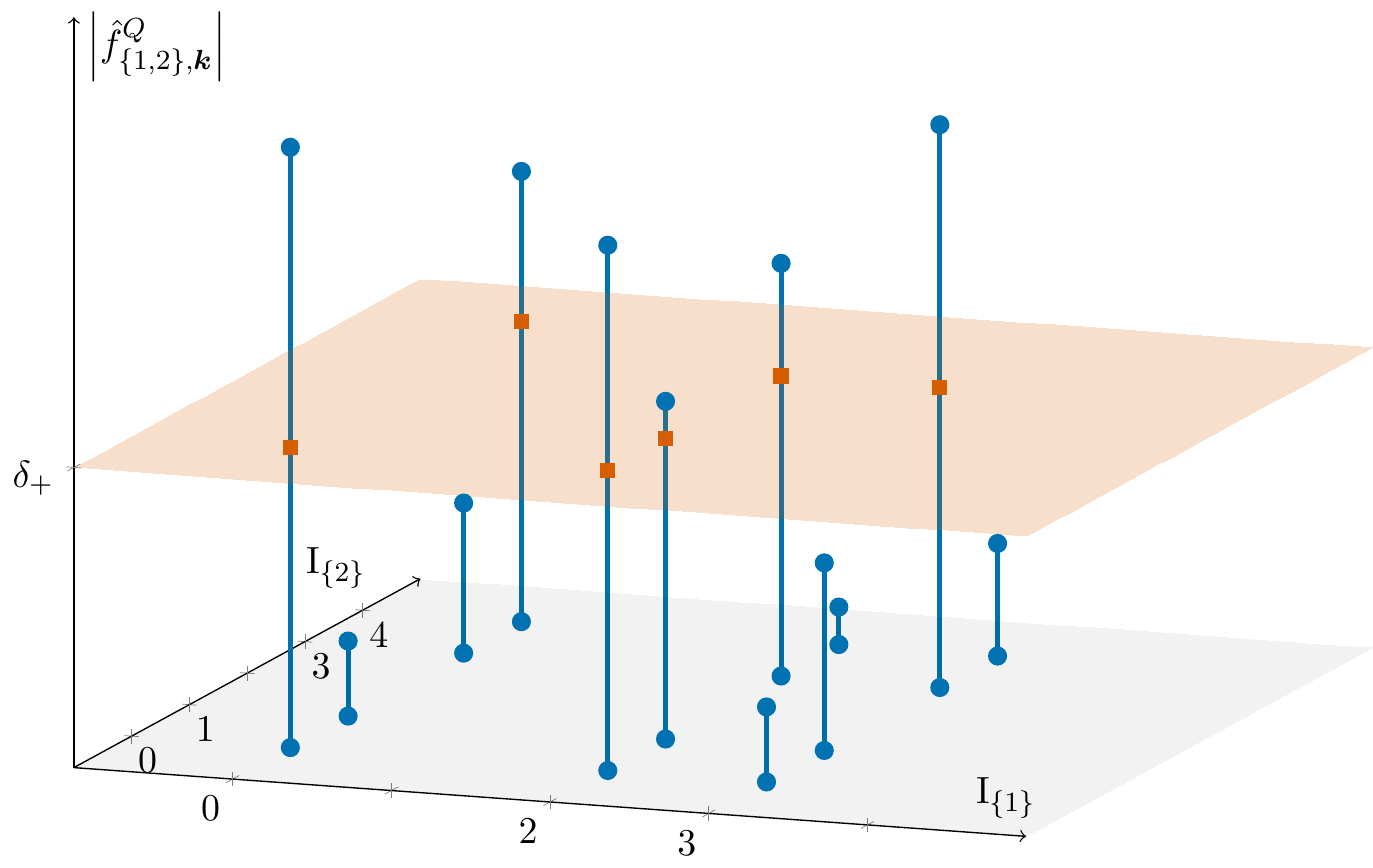}		
	\end{center}
	\caption{Two-dimensional detection step: Only the six indices $\boldk \in \N^2$ where $\fhat_{\lbrace 1,2\rbrace, \boldk}^Q \geq \delta_+$ holds are part of the new detected index set $\I_{\{1,2\}}$. All other two-dimensional indices $\boldk$ from the candidate set $K=\I_{\lbrace 1 \rbrace}\times\I_{\lbrace 2 \rbrace}$ are neglected after this step. For a clearer illustration, only one detection iteration $r=1$ is shown.} 
	\label{subfig:c}
\end{subfigure}
\newline
\begin{subfigure}[c]{0.95\textwidth}
	\begin{center}
		\includegraphics{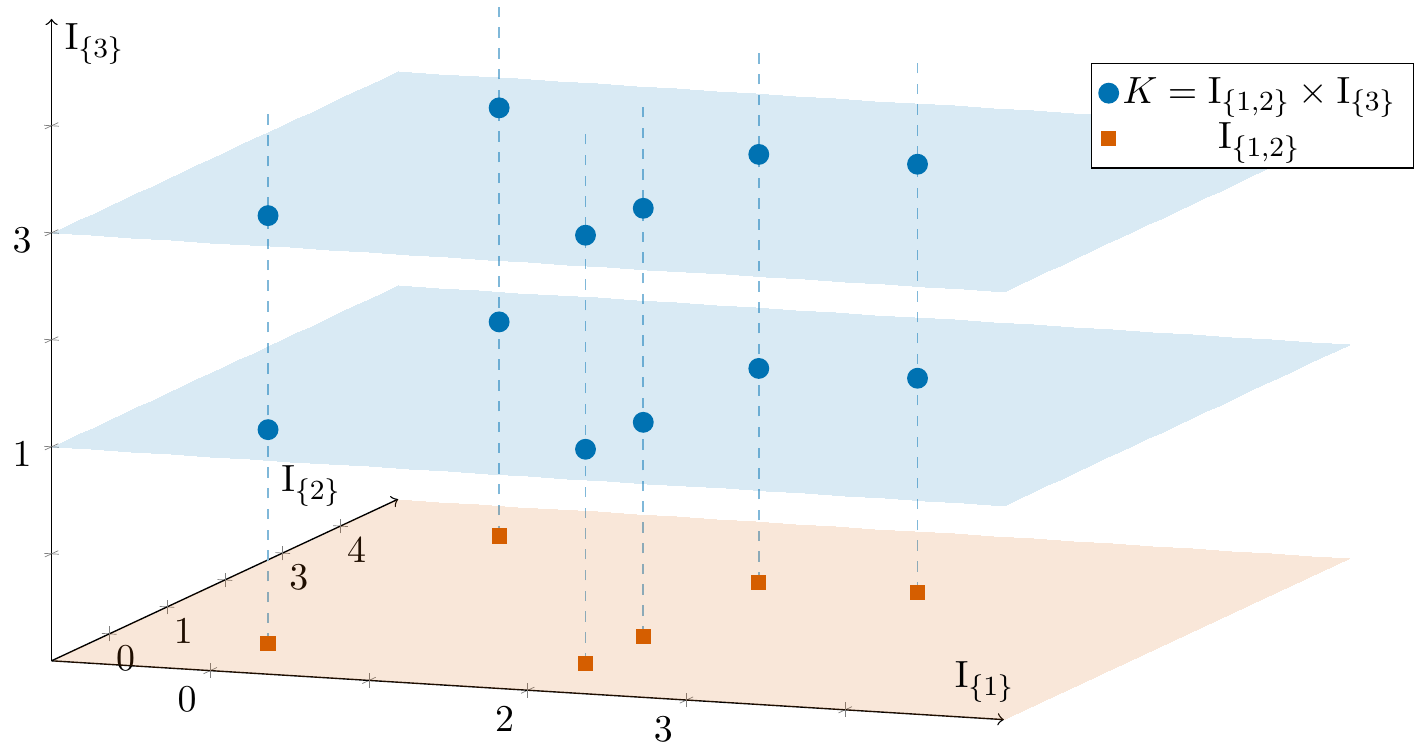}	
	\end{center}
	\caption{Three-dimensional candidate set: The new candidate set $K$ is now just the cartesian product of the two-dimensional detected set $\I_{\lbrace 1,2 \rbrace} = \{(0,0), (2,0), (2,1), (2,3), (3,3), (0,4)\}$ and the one-dimensional detected set $\I_{\lbrace 3 \rbrace}=\{1,3\}$.} 
	\label{subfig:d}
\end{subfigure}
\caption{Visualization of the first steps of Algorithm \ref{alg:main} for some example function $f$ using $\Gamma = \{0,1,2,3,4\}^d$.}
\label{fig:algo}
\end{figure}

\subsection{Complexity}\label{subsec:comp}

The sampling complexity as well as the computational complexity of Algorithm \ref{alg:main} obviously depends strongly on the function space $L_2(\D,\mu)$ as well as the used reconstruction methods in the corresponding function spaces $L_2(\D_{\u},\mu_{\u})$. To simplify the brief consideration of the complexity in this section, we assume that all function spaces $L_2(\D_j,\mu_j), j=1,\ldots,d,$ are equal and we can apply cubature method $Q$ for each product of those spaces.

In each part of Algorithm \ref{alg:main} the cubature method $Q$, see \eqref{eq:proj_coef_approx}, chooses an amount of sampling nodes, which depends on the amount of candidates $\abs{K}$, the current dimensionality $\abs{\u} \in \lbrace 1,\ldots,d \rbrace$ and the cubature method $Q$ itself. Hence, we will denote this sampling amount by $S_Q(\abs{K},\abs{\u})$. Note that this adds the implicit assumption that $Q$ acts independently of the structure of the candidate set $K$. 

In step 1, we have the dimensionality $\abs{\u}=\abs{\lbrace t \rbrace}=1$ and the candidate set $K = \P_{\u}(\Gamma) = \P_{\lbrace t \rbrace}(\Gamma)$. For most common choices of $\Gamma \in \N^d$, cf. Section \ref{subsec:gamma}, the one-dimensional projections $\P_{\lbrace t \rbrace}(\Gamma)$ are just the sets $\lbrace 0,\ldots,N-1\rbrace$ for some extension $N-1$, so $\abs{K}=N$. Since we sample $r$ times with different $\tilde{\boldx}$ in each dimension $t = 1,\ldots,d$, the amount of sampling nodes in step 1 of Algorithm \ref{alg:main} is then $drS_Q(N,1)$.

In step 2, we sample $r$ times for each dimensionality $\abs{\u} = 2,\ldots,d-1$ and $1$ time for $\abs{\u}=d$. The size of the $t$-th candidate set $\I_{\{1,\ldots,t-1\}} \times \I_{\{t\}}$ for $t=2,\ldots,d$ is bounded by $r^2 s_\mathrm{local}^2$, since both index sets contain at most $r s_\mathrm{local}$ indices (if the detected sets for each detection iteration were pairwise disjoint). The intersection with the projection of $\Gamma$ may only decrease the true number of samples. Hence, we end up with at most $S_Q(r^2s_\mathrm{local}^2,d) + \sum_{t=2}^{d-1}rS_Q(r^2s_\mathrm{local}^2,t)$ samples for step 2.

Finally, if $s \sim s_\mathrm{local}$, the sampling complexity of Algorithm \ref{alg:main} is then $$\OO{drS_Q(N,1) + \sum_{t=2}^{d-1}rS_Q(r^2s^2,t) + S_Q(r^2s^2,d)},$$ which is bounded by $\OO{drS_Q(r^2s^2,d)}$ if $r^2s^2 \geq N$ and $S_Q(\ast,d)\geq S_Q(\ast,t)$ hold for each $t=1,\ldots,d$.

For the computational complexity, we assume that all other steps of Algorithm \ref{alg:main} like the sampling itself, the choice of the random anchors $\tilde{\boldx}$ or the construction of the candidate sets $K$ are negligible. If we denote the computational complexity for the simultaneous numerical integration of all $\hat{f}_{\{u\},\boldk}^Q$ using the cubature method $Q$ by some expression $T_Q(\abs{K},\abs{\u})$ and assume no dependency on the structure of $K$ again, we receive the similar expression $$\OO{drT_Q(N,1) + \sum_{t=2}^{d-1}rT_Q(r^2s^2,t) + T_Q(r^2s^2,d)},$$ or with similar assumptions as before $\OO{drT_Q(r^2s^2,d)}$.

\subsection{A priori information}\label{subsec:gamma}

As already stated several times, we need $\Gamma$ to be large enough such that the desired indices we want to detect are all contained in it. Still, Algorithm \ref{alg:main} may benefits from  a better choice of $\Gamma$, so additional a priori information about the function $f$, since the amount of candidates especially in the higher-dimensional steps can be reduced significantly. While a full grid approach 
\begin{align*}
\Gamma = \hat{G}^d_n \coloneqq \left\lbrace \boldk \in \N^d: \vert\vert \boldk \vert\vert_{\infty} \leq n \right\rbrace
\end{align*}
with large enough $n$ will always work, for smoother functions $f$ with rapidly decaying basis coefficients $c_{\boldk}$ a weighted hyperbolic cross approach
\begin{align*}
\Gamma = H_{n}^{d,\boldgamma} \coloneqq \left\lbrace \boldk \in \N^d: \prod_{j=1}^d \max \left(1,\frac{k_j}{\gamma_j}\right) \leq n \right\rbrace
\end{align*}
with weight $\boldgamma = (\gamma_j)_{j=1}^d \in (0,\infty)^d$ is preferable. But even if the decay of the coefficients is relatively slow, an $\ell_p$ ball approach
\begin{align*}
\Gamma = B_{p,n}^{d,\boldgamma} \coloneqq \left\lbrace \boldk \in \N^d: \left(\sum_{j=1}^d \left(\frac{k_j}{\gamma_j}\right)^p\right)^{\frac1p} \leq n \right\rbrace
\end{align*}
with weight $\boldgamma = (\gamma_j)_{j=1}^d \in (0,\infty)^d$ can also reduce the amount of samples and computation time needed. For $p=\infty$, $\Gamma$ is the (weighted) full grid.

Another reasonable choice for practical examples comes from the sparsity-of-effects principle, which states that a system is usually dominated by main effects and low-order interactions. In our case, the principle means that the indices $\boldk$ with a rather small number of non-zero components $\vert\vert \boldk \vert\vert_0$ belong to the largest basis coefficients $c_\boldk$, as we already noticed when working with parametric PDEs in \cite{KaPoTa22}. This is also one of the main principles behind various low-order methods like the popular ANOVA decomposition, cf. \cite{KuSloWaWo10,PoSchmi19,schmischkediss} and the references therein, or the SHRIMP method, cf. \cite{XieShiSchWa22}. For such a case, a low-order approach
\begin{align*}
\Gamma_{\tilde{d}} \coloneqq \Gamma \cap \left\lbrace \boldk \in \N^d: \vert\vert \boldk \vert\vert_0 \leq \tilde{d} \right\rbrace
\end{align*}
with superposition dimension $\tilde{d} \in \N$ should be combined with any of the previous choices.

\begin{figure}
\begin{subfigure}[c]{0.32\textwidth}
	\begin{center}
		\includegraphics{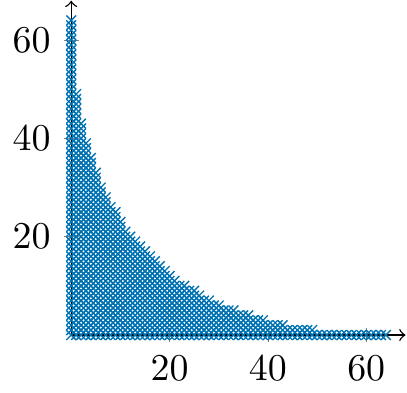}	
	\end{center}
	\caption{The $\ell_\frac12$ ball $B_{\frac12,64}^{2,\boldone}$ with weights $\boldone = (1,1)^\top$.}
	\label{subfig:gamma_b}
\end{subfigure}
~
\begin{subfigure}[c]{0.32\textwidth}
	\begin{center}
		\includegraphics{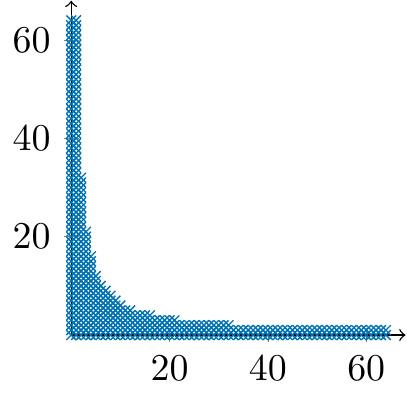}	
	\end{center}
	\caption{The hyperbolic cross $H_{64}^{2,\boldone}$ with weights $\boldone = (1,1)^\top$.} 
	\label{subfig:gamma_c}
\end{subfigure}
~
\begin{subfigure}[c]{0.32\textwidth}
	\begin{center}
		\includegraphics{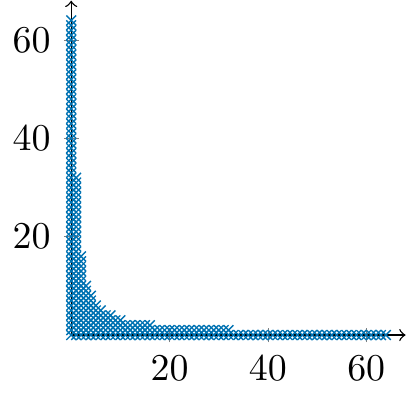}		
	\end{center}
	\caption{The hyperbolic cross $H_{128}^{2,\frac12\boldone}$ with weights $\frac12\boldone = (\frac12,\frac12)^\top$.} 
	\label{subfig:gamma_d}
\end{subfigure}
\caption{Illustration of three different choices for the search space $\Gamma$ in $d=2$ dimensions with similar extension.}
\label{fig:gamma}
\end{figure}

Table \ref{tbl:index} shows the size of the search space $\Gamma$ in $d=10$ dimensions for some examples with weights $\boldgamma = \boldone = (1,\ldots,1)^\top$ and their reduced versions $\Gamma_{\tilde{d}}$. Figure \ref{fig:gamma} illustrates three different index sets in two dimensions.

\begin{table}[tb]
\centering
\begin{tabular}{c|c||c|c|c} 
 \multicolumn{2}{c||}{$\tilde{d}$} & \multirow{2}{*}{none} & \multirow{2}{*}{$4$} & \multirow{2}{*}{$2$} \\
 $\Gamma$ & $n$ & & & \\
 \hline\hline
 \multirow{3}{*}{$\hat{G}^d_n$} & 8 & $3.49\cdot 10^9$ & $9.25\cdot 10^5$ & $2961$ \\\cline{2-5}
 & 16 & $2.02\cdot 10^{12}$ & $1.43\cdot 10^7$ & $1.17\cdot 10^4$ \\\cline{2-5}
 & 32 & $1.53\cdot 10^{15}$ & $2.24\cdot 10^8$ & $4.64\cdot 10^4$ \\
 \hline
 \multirow{3}{*}{$B_{2,n}^{d,\boldone}$} & $4$ & $4.32\cdot 10^4$ & $9161$ & $401$ \\\cline{2-5}
 & $6$ & $1.05\cdot 10^6$ & $5.64\cdot 10^4$ & $1051$ \\\cline{2-5}
 & $8$ & $1.19\cdot 10^7$ & $1.94\cdot 10^5$ & $1926$ \\
 \hline
 \multirow{3}{*}{$B_{1,n}^{d,\boldone}$} & $8$ & $4.38\cdot 10^4$ & $2.28\cdot 10^4$ & $1341$ \\\cline{2-5}
 & $12$ & $6.47\cdot 10^5$ & $1.33\cdot 10^5$ & $3091$ \\\cline{2-5}
 & $16$ & $5.31\cdot 10^6$ & $4.55\cdot 10^5$ & $5561$ \\
 \hline
 \multirow{3}{*}{$B_{\frac12,n}^{d,\boldone}$} & $32$ & $5.11\cdot 10^4$ & $5.00\cdot 10^4$ & $6891$ \\\cline{2-5}
 & $48$ & $3.44\cdot 10^5$ & $2.72\cdot 10^5$ & $1.61\cdot 10^4$ \\\cline{2-5}
 & $64$ & $1.55\cdot 10^6$ & $9.30\cdot 10^5$ & $2.90\cdot 10^4$ \\
 \hline
 \multirow{3}{*}{$H_{n}^{d,\boldone}$} & $8$ & $1.10\cdot 10^5$ & $1.88\cdot 10^4$ & $981$ \\\cline{2-5}
 & $16$ & $4.18\cdot 10^5$ & $5.85\cdot 10^4$ & $2411$ \\\cline{2-5}
 & $32$ & $1.52\cdot 10^6$ & $1.72\cdot 10^5$ & $5676$ \\
 \hline
\end{tabular}
\caption{Sizes of different search spaces $\Gamma$ for $d=10$ and their smaller versions with superposition dimensions $\tilde{d}=2$ and $\tilde{d}=4$.}
\label{tbl:index}
\end{table}

\subsection{Alternative increment strategies}\label{subsec:alternative}

One main feature of the dimension-incremental method is the combination of the detected, one-dimensional index set projections $\I_{\{t\}}, t=1,\ldots,d$. Algorithm \ref{alg:main}
realizes this in the most intuitive way by adding the dimension $t$ to the already detected set $\I_{\{1,\ldots,t-1\}}$ in each dimension increment $t=2,\ldots,d$. This classical approach, which we will call \textit{one-by-one strategy}, is sketched in Figure \ref{fig:obo} for $d=9$. The same approach was exploited in \cite{PoVo14}, where so-called reconstructing rank-1 lattices were used for the computation of the projected coefficients $\fhat_{\u,\boldk}^Q$. Therein, these lattices were computed component-by-component-wise and hence perfectly matched this incremental strategy. Another advantage of the one-by-one strategy can be seen when we combine the for-loop over $t$ in step 1 of Algorithm \ref{alg:main} with the for-loop over $t$ in step 2. This is possible, since each one-dimensional projection set $\I_{\{t\}}$ is only used in the corresponding $t$-th iteration in step 2. Hence, we only need one such one-dimensional projection set $\I_{\{t\}}$ at a time and may save additional memory space by simply overwriting the previous one in the next iteration $t+1$.

Obviously, there is no need to limit ourselves to this straight forward strategy in general. In the remaining part of this section, we therefore discuss some alternative increment strategies as well as possible advantages and disadvantages of these approaches. Note that all those strategies just yield several improvements, but are not ``optimal'' in any sense. Even worse, such ``optimalities'' heavily depends on the given problem and the corresponding dimension $d$, the cubature rules $Q$ and the used algorithm parameters. Hence, it is very tricky to come up with an overall good strategy for every possible setting. Also, the reader needs to decide on its own, which kind of ``optimality'' is even aimed for, e.g., small complexities of the cubature rules, memory efficiency or possible parallelizations.

\subsubsection*{Dyadic strategy}

In step 2 of Algorithm \ref{alg:main}, we most importantly construct the cubature rule $Q$ and compute the projected coefficients $\fhat_{\u,\boldk}^Q$. Since the values for $\tilde{\boldx}$ are random but fixed, these are basically $\abs{\u}$-dimensional problems. Depending on the used cubature rule, the size $\abs{\u}$ might heavily influences the amount of cubature nodes $M$ as well as the computational costs.

\begin{figure}[tb]
	\centering	
	\begin{subfigure}[c]{0.95\textwidth}
	\centering
	\resizebox{0.95\textwidth}{!}{
	\includegraphics{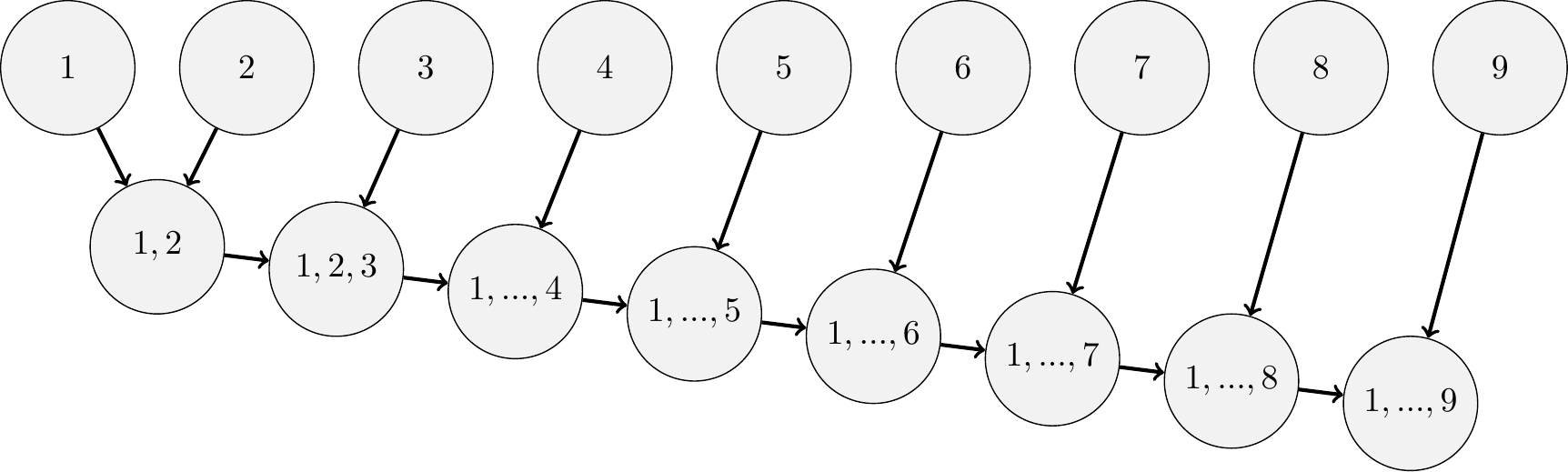}	
	}
	\subcaption{One-by-one strategy}\label{fig:obo}
	\end{subfigure}
	\par\medskip	
	\begin{subfigure}[c]{0.95\textwidth}
	\centering
	\resizebox{0.95\textwidth}{!}{
	\includegraphics{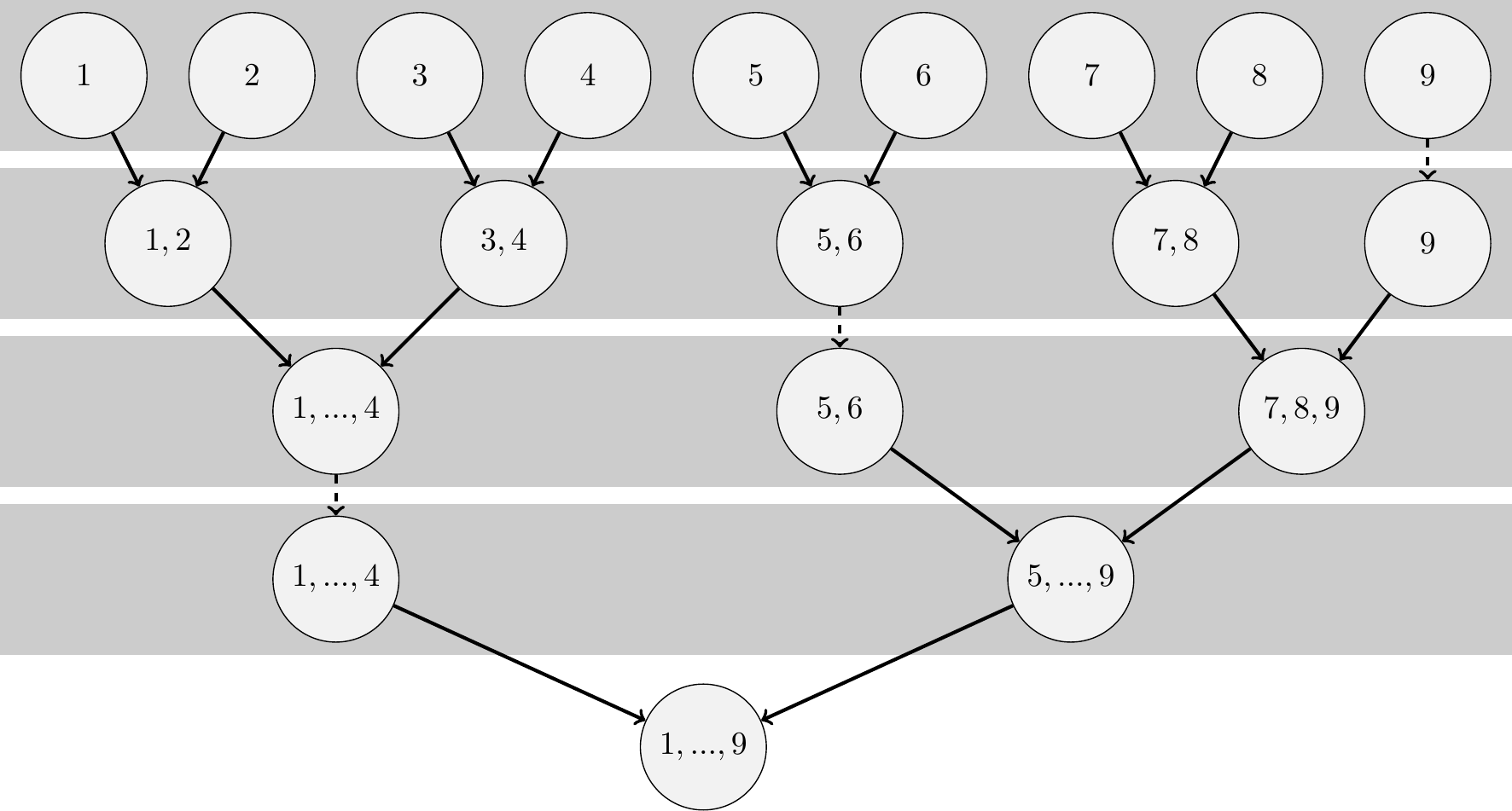}	
	}
	\subcaption{Dyadic strategy}\label{fig:dad}
	\end{subfigure}
	\caption{Visualization of non-data-driven increment strategies for $d=9$.}
  \label{fig:increment_strategies}
\end{figure}

Using the one-by-one strategy, we consider $d-1$ dimension-incremental steps, where each dimensionality $\abs{\u} = t \in \{2,\ldots,d\}$ appears exactly once. The \textit{dyadic strategy} now aims for more incremental steps with lower dimensionality while keeping the overall number of steps constant. This strategy combines the two projected index sets $\I_\u$ and $\I_\v$ with smallest dimensionalities $\abs{\u}$ and $\abs{\v}$ in each step. If there are several sets of same dimensionalities, e.g., at the beginning when there are $d$ sets with dimension $\abs{\u} = 1$, the set is randomly chosen among these candidates. Rearranging these dimension-incremental steps into stages as in Figure \ref{fig:dad}, the dyadic structure can be seen. Note that for $d\not=2^k, k\in\N$, some stages have to keep one projected index set untouched since there was an odd number of sets to combine, which will be the projected index set $\I_\u$ with the highest dimensionality $\abs{\u}$. This is the case for $\I_{\{9\}}$ in the first stage, $\I_{\{5,6\}}$ in the second stage and $\I_{\{1,\ldots,4\}}$ in the third stage in Figure \ref{fig:dad}, visualized using the dashed arrow.

As mentioned above, this strategy reduces the dimensionalities in many steps tremendously for large $d$. Even for the relatively small $d=9$ in Figure \ref{fig:dad}, the dyadic strategy uses four $2$-dimensional steps and only one $3$-, $4$-, $5$- and $9$-dimensional step each instead of one $t$-dimensional step for each $t=2,\ldots,9$ as in the one-by-one strategy. The additional computational effort for the realization of this strategy is also relatively small, since it depends only on the particular $d$ used and can even be precomputed. Finally, many dimension-incremental steps can be performed in parallel, as they are not dependent on each other, which allows additional time savings. Unfortunately, the dyadic strategy is not implementable as memory-friendly as the one-by-one strategy since some of the projected index sets $\I_\u$ block memory for several steps while they are not picked for the next combination.

\subsubsection*{Data-driven one-by-one strategy}

While the dyadic strategy aims for smaller dimensionalities $\abs{\u}$ in the dimension-incremental steps, the size of the candidate set $K \subset \N^{\abs{\u}}$, so how many approximated projected coefficients $\fhat_{\u,\boldk}^Q$ we need to compute, is probably another crucial factor influencing the performance of the cubature rules $Q$ and hence the overall performance of our algorithm. In general, we have $K = \left(\I_\u \times \I_\v\right) \cap \P_{\u\cup\v}(\Gamma)$, e.g., $\I_\u = \I_{\{1,\ldots,t\}}$ and $\I_\v = \I_{\{t\}}$ in a dimension-increment in Algorithm \ref{alg:main}. The size $\abs{K}$ depends mainly on the size of the projected index sets $\I_\u$ and $\I_\v$ it is built from. Therefore, we now present two so-called \textit{data-driven} strategies, where these sizes $\abs{\I_\u}$ are examined before each dimension-incremental step and then it is decided, which $\I_\u$ and $\I_\v$ to choose for the next step. Note that an investigation of the sizes of all possible $K$ instead of just the sizes of all available $\I_\u$ might be even more favorable, especially for challenging choices of $\Gamma$, but also needs even more computational effort and is therefore not considered herein.

The following approach is based on the classical one-by-one strategy and is therefore called \textit{data-driven one-by-one strategy}. We start with the computation of all the one-dimensional projected index sets $\I_{\{t\}}, t=1,\ldots,d$. We proceed as in the one-by-one strategy, but instead of working through the sets lexicographically, so from $\I_{\{1\}}$ to $\I_{\{d\}}$, we rearrange them, ordered by their descending size. In particular, we go from $\I_{\{t_1\}}$ to $\I_{\{t_d\}}$ with $t_i \in \{1,\ldots,d\}$ for all $i \in \{1,\ldots,d\}$ and $\vert\I_{\{t_i\}}\vert \geq \vert\I_{\{t_j\}}\vert, 1\leq j \leq i \leq d$. For instance, in Figure \ref{fig:dobo} the set $\I_{\{5\}}$ is the largest one and therefore the starting point of the data-driven one-by-one strategy.

The advantage of this approach compared to the classical one-by-one strategy is the fact that the higher-dimensional candidate sets $K$ in the later dimension-incremental steps are probably smaller and hence the construction of the cubature rule $Q$ as well as the computation of the coefficients might be considerably cheaper in terms of sampling points and computation time.

\subsubsection*{Data-driven dyadic strategy}

\begin{figure}[tb]
	\centering	
	\begin{subfigure}[c]{0.95\textwidth}
	\centering
	\resizebox{0.95\textwidth}{!}{
	\includegraphics{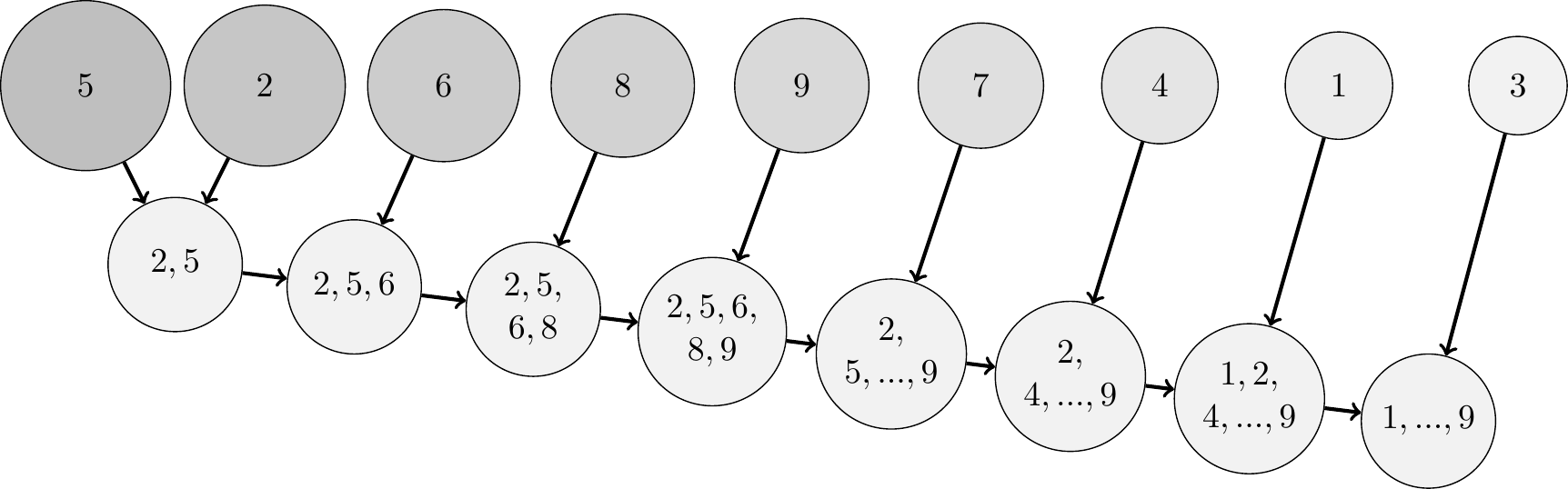}	
	}
	\subcaption{Data-driven one-by-one strategy}\label{fig:dobo}
	\end{subfigure}
	\par\medskip
	\begin{subfigure}[c]{0.95\textwidth}
	\centering
	\resizebox{0.95\textwidth}{!}{
	\includegraphics{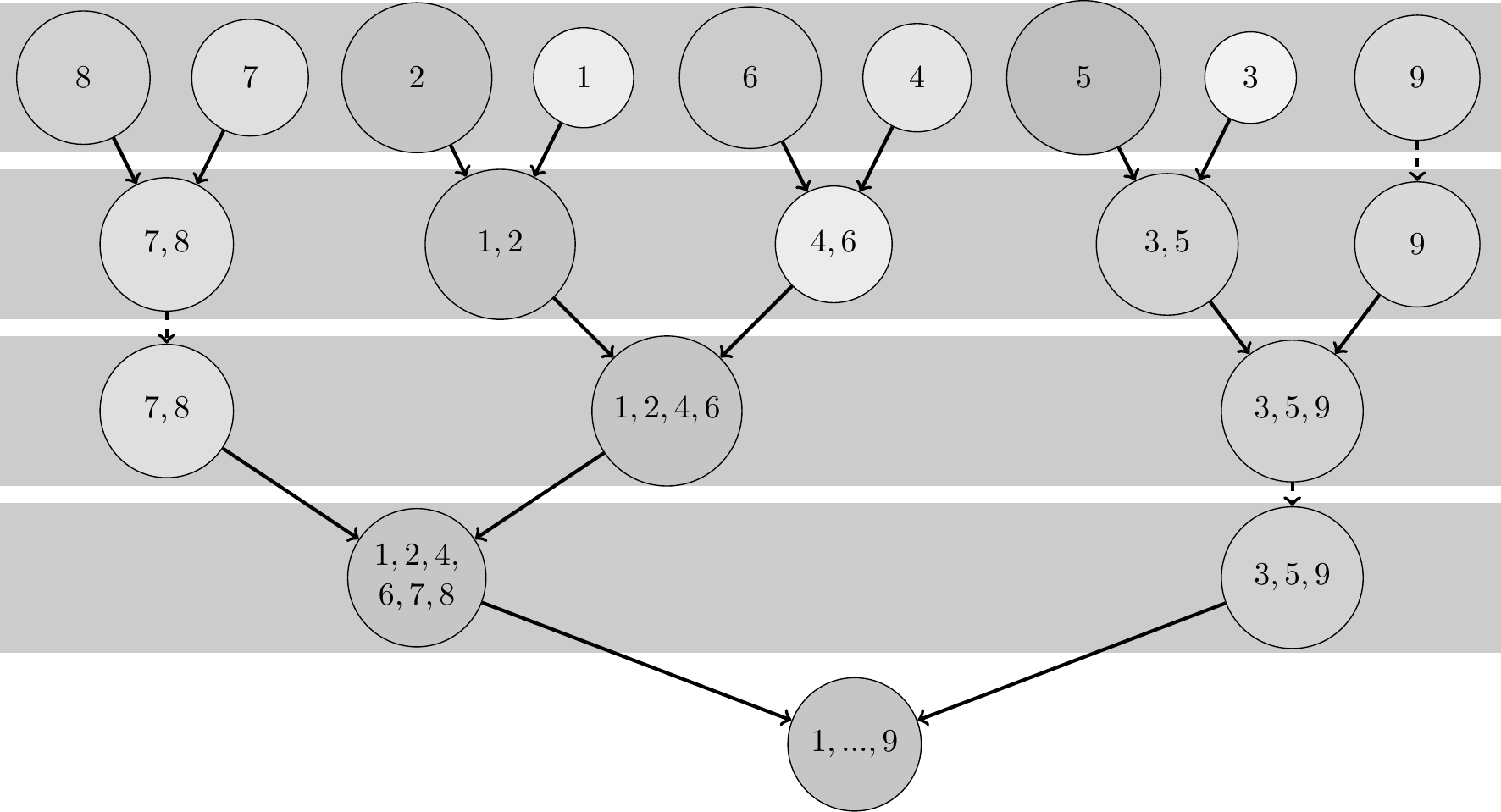}	
	}
	\subcaption{Data-driven dyadic strategy}\label{fig:ddd}
	\end{subfigure}
	\caption{Visualization of data-driven increment strategies for $d=9$.}
  \label{fig:increment_strategies_dd}
\end{figure}

Finally, the \textit{data-driven dyadic strategy} combines the advantages of both the dyadic strategy and the data-driven one-by-one strategy, illustrated for $d=9$ in Figure \ref{fig:ddd}. While the rearrangement into stages in the dyadic strategy was just for visualization, it is an essential part of the strategy now. In the first stage, we investigate the size of all one-dimensional sets $\I_{\{t\}}, t=1,\ldots,d$, and perform $\floor{d/2}$ dimension-incremental steps, where we combine the largest set with the smallest one, the second largest with the second smallest and so on. In the next stage, we then take a look at all those new, two-dimensional index sets (and the possibly leftover one-dimensional set from the first stage if $d$ is odd) and perform multiple dimension-incremental steps again, using the same criterion as in the first stage. This strategy is then repeated as often as needed, so until we end up with the full index set $\I_{\{1,\ldots,d\}}$. Note that the dimensionality $\abs{\u}$ of the sets is not considered in this kind of strategy, but higher-dimensional sets are more likely to be larger as well and is therefore involved implicitly. If the number of available sets is odd in any stage, the median sized set is just kept as it is for the next stage, e.g., sets $\I_{\{9\}}$, $\I_{\{7,8\}}$ and $\I_{\{3,5,9\}}$ in the first three stages in Figure \ref{fig:ddd}.

As for the dyadic strategy, we manage to end up with less high-dimensional combination steps as for the one-by-one or data-driven one-by-one strategy, which can be performed in parallel again. On the other hand, the additional size criterion in each stage avoids cases, where two very large sets are combined and hence the size of the corresponding candidate set $K$ grows unnecessarily large. This can easily happen in the dyadic strategy since it only considers the dimensionalities $\abs{\u}$ of the sets $\I_{\u}$ and also chooses randomly among sets with the same dimensionality.

\section{Theoretical detection guarantee}\label{sec:theory}

In this section, we show a bound on the number of detection iterations $r$ in Algorithm \ref{alg:main} such that we can ensure the successful detection of all indices $\boldk$ belonging to basis coefficients $c_{\boldk}$, whose magnitude is larger than an absolute threshold $\delta$, with high probability. Therefore, we follow the main steps in \cite{KaPoVo17} and generalize their theoretical results. As explained in Section \ref{subsec:setup}, we consider smooth enough, multivariate functions $f \in L_2(\D,\mu)$ of the form
\begin{align*}%\label{eq:function}
f(\boldx) = \sum_{\boldk \in \N^d} c_{\boldk}\Phi_{\boldk}(\boldx)
\end{align*}
for some coefficients $\lbrace c_{\boldk}\rbrace_{\boldk\in\N^d}$. Further, we denote with $\I_{\delta} \subset \N^d$ the unknown, finite index set such that
\begin{align}\label{eq:delta_condition}
\begin{cases} \abs{c_{\boldk}} \geq \delta & \forall \boldk \in \I_{\delta} \\
\abs{c_{\boldk}} < \delta & \forall \boldk \not\in \I_{\delta} \end{cases}
\end{align}
holds and $\Gamma \subset \N^d$ as before the suitable search domain containing this index set $\I_{\delta}$, i.e., $\I_\delta \subset \Gamma$.

We compute the adaptive approximation in $d$ dimension increment steps in Algorithm \ref{alg:main}. In each of the dimension increment steps at most three probabilistic sub-steps are performed: \begin{itemize}
\item The detection of the one-dimensional projections $\I_{\{t\}}$ in step 1, which is successful, if the event $E_{1,t} \coloneqq \left\lbrace \mathcal{P}_{\u}(\I_\delta) \subset \I_{\u} \right\rbrace$ with $\u = \{t\}$ occurs.
\item The (possibly probabilistic) construction of the cubature rule $Q$ for some index set $K$ in step 2, i.e., the successful computation of $M, (w_j)_{j=1}^M$ and $(\boldxi_j)_{j=1}^M$. We define $$E_{2,t} \coloneqq \left\lbrace \text{Successful construction of $Q$ for $K$} \right\rbrace.$$
\item The detection of the multi-dimensional projections $\I_{\{1,\ldots,t\}}$ in step 2, which is successful, if the event $E_{3,t} \coloneqq \left\lbrace \mathcal{P}_{\u} (\I_\delta) \subset \I_{\u} \right\rbrace$ with $\u = \{1,\ldots,t\}$ occurs. 
\end{itemize} 

If each of these probabilistic sub-steps is successful, we detect all indices from $\I_\delta$. We use the union bound to estimate the corresponding probability
\begin{align}\label{eq:prob_main}
\P \left( \bigcap_{t=1}^d E_{1,t} \cap \bigcap_{t=2}^d E_{2,t} \cap \bigcap_{t=2}^d E_{3,t} \right) &= 1 - \P \left( \bigcup_{t=1}^d E_{1,t}^\comp \cup \bigcup_{t=2}^d E_{2,t}^\comp \cup \bigcup_{t=2}^d E_{3,t}^\comp \right) \nonumber \\
&\geq 1 - \sum_{t=1}^d \P(E_{1,t}^\comp) - \sum_{t=2}^d \P(E_{2,t}^\comp) - \sum_{t=2}^d \P(E_{3,t}^\comp).
\end{align}
We aim for a failure probability $\varepsilon \in (0,1)$ of the whole algorithm. We split this up such that each probabilistic sub-step has an equal upper bound on its failure probability of $\varepsilon/(3d)$. Hence, we now estimate the probabilities $\P(E_{1,t}^\comp)$ and $\P(E_{3,t}^\comp)$. Upper bounds on $\P(E_{2,t}^\comp)$ depend on the used cubature rule $Q$ and its construction and are therefore not considered here.

\subsection{Failure probability estimates}

First, we recall the definition and estimate of the approximated projected coefficients \eqref{eq:proj_coef_approx} as well as the projection error term \eqref{eq:noise} from Section \ref{subsec:pc} and apply them for $J = \I_\delta$. We use \begin{itemize}
\item $\u=\{t\}$ and $K = \mathcal{P}_{\u}(\Gamma) \subset \N$ for the one-dimensional projections in step 1 and
\item $\u = \{1,\ldots,t\}$ and $K =(\I^{(1,\ldots,t-1)} \times \I^{(t)})\cap\mathcal{P}_{\u}(\Gamma) \subset \N^t$ for the multi-dimensional projections in step 2 
\end{itemize}
of Algorithm \ref{alg:main}. In each case, we get for $\boldk \in K$ the formula
\begin{align}\label{eq:proj_coef_theory}
\fhat_{\u,\boldk}^Q (\tilde{\boldx}) &= \sum_{\boldh = (\boldk,\boldh_{\u^\comp})_\u \in J} c_{\boldh} \Phi_{\u^\comp,\boldh_{\u^\comp}}(\tilde{\boldx}) + \Psi_{\u,\boldk}^{Q,J}(\tilde{\boldx}).
\end{align}
Note that the first part of \eqref{eq:proj_coef_theory} is independent of the used cubature rule $Q$. Hence, the approximated projected coefficients $\fhat_{\u,\boldk}^Q$ are also independent of $Q$ up to the projection error term $\Psi_{\u,\boldk}^{Q,J}$, depending on the basis coefficients of $f$ in $\I_\delta^\comp$, so all coefficients $c_\boldk$ with absolute value smaller than $\delta$. 

The following Lemma based on \cite[Lem.\ 4]{KaPoVo17} gives an estimate on the probability that such sums \eqref{eq:proj_coef_theory} take small function values for randomly drawn anchors $\tilde{\boldx}$.

\begin{Lemma}\label{lem:trig_prob}
Consider $\v \subset \lbrace 1,\ldots,d \rbrace$ and the space $\D_\v \coloneqq \bigtimes_{j\in\v} \D_j \subset \R^{\abs{\v}}$ with the corresponding product measure $\mu_\v$ and basis functions $\Phi_{\v,\boldh}$, which are bounded by the finite constant $B_\v$, cf.\ Section \ref{subsec:pc}. 

Let a function $g: \D_\v \rightarrow \C,\, g(\boldx) \coloneqq \sum_{\boldh\in\tilde{\I}} \hat{g}_{\boldh} \Phi_{\v,\boldh}(\boldx) \not\equiv 0,\, \tilde{\I} \subset \N^{\abs{\v}}, \vert\tilde{\I}\vert < \infty,$ be given and let $\Psi: \D_\v \rightarrow \C$ be some function with $\norm{\Psi}_{L_1(\D_\v)} < \delta_\Psi$. Moreover, let $X_j\in \D_j,\, j \in \v,$ be independent, $\mu_j$-distributed random variables and we denote by $\boldX \coloneqq (X_j)_{j\in \v} \in \D_\v$ the random vector. 

If $\max_{\boldh\in\tilde{\I}}\vert\hat{g}_{\boldh}\vert > B_\v(\delta_+ + \delta_\Psi)$ for some $\delta_+ >0$, then
\begin{align}\label{eq:l1_zwei}
\P(\abs{(g+\Psi)(\boldX)}<\delta_+)\leq \left( 1 - \frac{B_\v^{-1} \max_{\boldh \in \tilde{\I}} \abs{\ghat_\boldh}-\delta_\Psi-\delta_+}{B_\v \sum_{\boldh \in \tilde{\I}} \abs{\ghat_\boldh}+\norm{\Psi}_{L_\infty(\D_\v)}} \right) \eqqcolon q < 1.
\end{align}
Choosing $r$ random vectors $\boldX_1,\ldots,\boldX_r \in\D_\v$ independently, we observe
\begin{align}\label{eq:l1_drei}
\P\left(\bigcap_{i=1}^r \left\lbrace\abs{(g+\Psi)(\boldX_i)}<\delta_+\right\rbrace\right) \leq q^r.
\end{align}
\end{Lemma}
\begin{proof}
We proceed as in \cite{KaPoVo17} and refer to the lower bound
\begin{align*}
\frac{\E h(Y) - h(\delta_+)}{\norm{h(Y)}_{L_\infty(\R)}} \leq \P(\abs{Y} \geq \delta_+)
\end{align*}
from \cite[Par.\ 9.3.A]{Lo77}. Applying this for the even and on $[0,\infty)$ nondecreasing function $h(t) \coloneqq \abs{t}$ and $Y=(g+\Psi)(\boldX)$ leads to
\begin{align*}
\P(\abs{(g+\Psi)(\boldX)}\geq\delta_+) &\geq \frac{\norm{g+\Psi}_{L_1(\D_\v)}-\delta_+}{\norm{g+\Psi}_{L_\infty(\D_\v)}} \\
&\geq \frac{\norm{g}_{L_1(\D_\v)}-\norm{\Psi}_{L_1(\D_\v)}-\delta_+}{\norm{g}_{L_\infty(\D_\v)}+\norm{\Psi}_{L_\infty(\D_\v)}} \\
&\geq \frac{\norm{g}_{L_1(\D_\v)}-\delta_\Psi-\delta_+}{\norm{g}_{L_\infty(\D_\v)}+\norm{\Psi}_{L_\infty(\D_\v)}}
\intertext{and consequently}
\P(\abs{(g+\Psi)(\boldX)}\leq\delta_+) &\leq 1-\frac{\norm{g}_{L_1(\D_\v)}-\delta_\Psi-\delta_+}{\norm{g}_{L_\infty(\D_\v)}+\norm{\Psi}_{L_\infty(\D_\v)}}.
\end{align*}
Since 
\begin{align*}
\abs{\hat{g}_{\boldh}} &= \abs{\int_{\D_\v} g(\boldx)\,\overline{\Phi_{\v,\boldh}(\boldx)}\diff\mu_\v(\boldx)} \\
&\leq \int_{\D_\v} \abs{g(\boldx)} \abs{\Phi_{\v,\boldh}(\boldx)} \diff\mu_\v(\boldx) \\
&\leq \norm{\Phi_{\v,\boldh}}_{L_\infty(\D_\v)} \int_{\D_\v} \abs{g(\boldx)} \diff\mu_\v(\boldx) = \norm{\Phi_{\v,\boldh}}_{L_\infty(\D_\v)} \norm{g}_{L_1(\D_\v)}
\end{align*}
for all $\boldh \in \tilde{\I}$, we have
\begin{align*}
\max_{\boldh\in\tilde{\I}} \abs{\hat{g}_{\boldh}} \leq \norm{g}_{L_1(\D_\v)} \max_{\boldh\in\tilde{\I}} \norm{\Phi_{\v,\boldh}}_{L_\infty(\D_\v)} = \norm{g}_{L_1(\D_\v)} B_\v. 
\end{align*}
Together with the estimate 
\begin{align*}
\norm{g}_{L_\infty(\D_\v)} = \esssup_{\boldx\in\D_\v} \abs{g(\boldx)} &\leq \esssup_{\boldx\in\D_\v} \sum_{\boldh\in\tilde{\I}} \abs{\hat{g}_{\boldh}} \abs{\Phi_{\v,\boldh}(\boldx)}\\
&\leq \sum_{\boldh\in\tilde{\I}} \abs{\hat{g}_{\boldh}} \esssup_{\boldx\in\D_\v}  \abs{\Phi_{\v,\boldh}(\boldx)} \\
&\leq \sum_{\boldh\in\tilde{\I}} \abs{\hat{g}_{\boldh}} \norm{\Phi_{\v,\boldh}}_{L_\infty(\D_\v)} \leq B_\v \sum_{\boldh\in\tilde{\I}} \abs{\hat{g}_{\boldh}},
\end{align*}
we conclude
\begin{align*}
\P(\abs{(g+\Psi)(\boldX)}\leq\delta_+) &\leq 1-\frac{\norm{g}_{L_1(\D_\v)}-\delta_\Psi-\delta_+}{\norm{g}_{L_\infty(\D_\v)}+\norm{\Psi}_{L_\infty(\D_\v)}} \\
&\leq 1 - \frac{B_\v^{-1} \max_{\boldh \in \tilde{\I}} \abs{\ghat_\boldh}-\delta_\Psi-\delta_+}{B_\v \sum_{\boldh \in \tilde{\I}} \abs{\ghat_\boldh}+\norm{\Psi}_{L_\infty(\D_\v)}}.
\end{align*}
Using the assumption $\max_{\boldh\in\tilde{\I}} \abs{\hat{g}_{\boldh}} > B_\v(\delta_+ + \delta_\Psi)$, the estimate \eqref{eq:l1_zwei} holds and \eqref{eq:l1_drei} then follows directly.
\end{proof}

In Algorithm \ref{alg:main}, we compute the approximated projected coefficients $\fhat_{\u,\boldk}^Q (\tilde{\boldx})$ for all $\boldk$ in our candidate sets $\mathcal{P}_\u(\Gamma)$ in step 1 and $\left( \I_{\{1,\ldots,t-1\}} \times \I_{\{t\}} \right) \cap \mathcal{P}_{\u}(\Gamma)$ in step 2.
Now, we apply Lemma \ref{lem:trig_prob} to those coefficients with $\boldk \in \mathcal{P}_\u(\I_\delta)$ and $\boldk \in \left( \I_{\{1,\ldots,t-1\}} \times \I_{\{t\}} \right) \cap \mathcal{P}_{\u}(\I_\delta)$, respectively, so those projected coefficients, we want to detect. This yields us bounds on the probability that they are below the detection threshold $\delta_+$ and therefore not detected by Algorithm \ref{alg:main}.

\begin{Corollary}\label{cor:single}
Let a threshold value $\delta > 0$ and a smooth enough function $f \in L_2(\D,\mu)$ be given. We consider the finite index set $\I_\delta$ such that \eqref{eq:delta_condition} holds.
\begin{itemize}
\item For fixed $1 \leq t \leq d$, $\u = \{t\}$, we denote $K(S) \coloneqq \mathcal{P}_\u(S), S \in \{\I_\delta,\Gamma\}$ and compute the one-dimensional approximated projected coefficients for the $t$-th component and
\item for fixed $1 < t < d$, $\u = \{1,\ldots,t\}$, we denote $K(S) \coloneqq \left( \I_{\{1,\ldots,t-1\}} \times \I_{\{t\}} \right) \cap \mathcal{P}_{\u}(S), S \in \{\I_\delta,\Gamma\}$ and compute the multi-dimensional approximated projected coefficients for the components $(1,\ldots,t)$
\end{itemize}
using a cubature rule $Q$ with \eqref{eq:quadrature_exact} and \eqref{eq:weight_sum} by
\begin{align*}
\fhat_{\u,\boldk}^Q (\tilde{\boldx}) = \sum_{j=1}^M w_j f(\boldxi_j,\tilde{\boldx})_\u\overline{\Phi_{\u,\boldk}(\boldxi_j)} \qquad \boldk \in K(\Gamma),
\end{align*}
where the anchors $\tilde{\boldx}_1,\ldots,\tilde{\boldx}_r \in \D_{\u^\comp}, r \in \N$ are independently chosen according to the corresponding product measure $\mu_{\u^\comp}$ at random. Further, we assume that there exists a $\delta_\Psi > 0$ such that $\norm{\Psi_{\u,\boldk}^{Q,I_\delta}}_{L_1(\D_{\u^\comp})} < \delta_\Psi$ holds for all $\boldk \in K(\I_\delta)$. 

Then, for $\delta_+>0$ with $B(\delta_+ + \delta_\Psi) < \min_{\boldk\in\I_\delta} \abs{c_\boldk}$ and $\boldk \in K(\I_\delta)$, the probability estimate
\begin{align}\label{eq:AtC}
&\P\left( \bigcup_{\boldk \in K(\I_\delta)} \left\lbrace \max_{\nu = 1,\ldots,r} \abs{\fhat_{\u,\boldk}^Q (\tilde{\boldx}_\nu)} < \delta_+ \right\rbrace \right) \leq \abs{\I_\delta} \left( \max_{\boldk \in K(\I_\delta)} q_{\u,\boldk} \right)^r
\end{align}
holds, where
\begin{align*}
q_{\u,\boldk} \coloneqq 1 - \frac{\displaystyle B^{-1} \max_{\boldh=(\boldk,\boldh_{\u^\comp})_\u\in\I_\delta} \abs{c_{\boldh}} - \delta_\Psi - \delta_+}{\displaystyle B \sum_{\boldh=(\boldk,\boldh_{\u^\comp})_\u\in\I_\delta} \abs{c_{\boldh}} + B^2 C_Q \sum_{\boldh \in \I_\delta^\comp} \abs{c_{\boldh}}} < 1,
\end{align*}
holds.

\begin{proof}
The estimate
\begin{align*}
\P \left( \abs{\fhat_{\u,\boldk}^Q (\tilde{\boldx})} < \delta_+ \right) &\leq 1 - \frac{\displaystyle B^{-1} \max_{\boldh=(\boldk,\boldh_{\u^\comp})_\u\in\I_\delta} \abs{c_{\boldh}} - \delta_\Psi - \delta_+}{\displaystyle B \sum_{\boldh=(\boldk,\boldh_{\u^\comp})_\u\in\I_\delta} \abs{c_{\boldh}} + \norm{\Psi_{\u,\boldk}^{Q,\I_\delta}}_{L_\infty(\D_{\u^\comp})}} \\
&\leq \left( 1 - \frac{\displaystyle B^{-1} \max_{\boldh=(\boldk,\boldh_{\u^\comp})_\u\in\I_\delta} \abs{c_{\boldh}} - \delta_\Psi - \delta_+}{\displaystyle B \sum_{\boldh=(\boldk,\boldh_{\u^\comp})_\u\in\I_\delta} \abs{c_{\boldh}} + B^2 C_Q \sum_{\boldh \in \I_\delta^\comp} \abs{c_{\boldh}}} \right) = q_{\u,\boldk} < 1
\end{align*}
or
\begin{align*}
\P \left( \abs{\fhat_{\u,\boldk}^Q (\tilde{\boldx})} < \delta_+ \right) &\leq q_{\u,\boldk} < 1
\end{align*}
holds due to \eqref{eq:l1_zwei} in Lemma \ref{lem:trig_prob} and the estimate
\begin{align*}
\norm{\Psi_{\u,\boldk}^{Q,\I_\delta}}_{L_\infty(\D_{\u^\comp})} &= \esssup_{\tilde{\boldx}\in\D_{\u^\comp}} \abs{\sum_{\boldh \in \I_\delta^\comp} c_{\boldh}\,\Phi_{\u^\comp,\boldh_{\u^\comp}}(\tilde{\boldx}) \sum_{j=1}^M w_j \Phi_{\u,\boldh_\u}(\boldxi_j) \,\overline{\Phi_{\u,\boldk}(\boldxi_j)}} \\
&\leq \sum_{\boldh \in \I_\delta^\comp} \abs{c_{\boldh}} \esssup_{\tilde{\boldx}\in\D_{\u^\comp}} \abs{\Phi_{\u^\comp,\boldh_{\u^\comp}}(\tilde{\boldx})} \sum_{j=1}^M \abs{w_j} \abs{\Phi_{\u,\boldh_\u}(\boldxi_j)}\abs{\Phi_{\u,\boldk}(\boldxi_j)} \\
&\leq B^2 C_Q \sum_{\boldh \in \I_\delta^\comp} \abs{c_{\boldh}}.
\end{align*}
Repeating the computation of $\fhat_{\u,\boldk}^Q (\tilde{\boldx})$ for independent randomly chosen anchors $\tilde{\boldx} = \tilde{\boldx}_1,\ldots,$ $\tilde{\boldx}_r \in \D_{\u^\comp}, r\in\N,$ we estimate
\begin{align*}
\P\left( \max_{\nu = 1,\ldots,r} \abs{\fhat_{\u,\boldk}^Q (\tilde{\boldx}_\nu)} < \delta_+ \right) \leq (q_{\u,\boldk})^r.
\end{align*}
Applying the union bound yields
\begin{align*}
&\P\left( \bigcup_{\boldk \in K(\I_\delta)} \left\lbrace \max_{\nu = 1,\ldots,r} \abs{\fhat_{\u,\boldk}^Q (\tilde{\boldx}_\nu)} < \delta_+ \right\rbrace \right) \leq \sum_{\boldk \in K(\I_\delta)} (q_{\u,\boldk})^r \leq \abs{\I_\delta} \left( \max_{\boldk \in K(\I_\delta)} q_{\u,\boldk} \right)^r.
\end{align*}
\end{proof}
\end{Corollary}

Finally, the following Lemma is based on \cite[Lem.\ 7]{KaPoVo17} and gives an estimate on the choice of the number of detection iterations $r$ such that the probabilities $\P(E_{1,t}^\comp)$ and $\P(E_{3,t}^\comp)$ in \eqref{eq:prob_main} are bounded by the desired value $\varepsilon/(3d)$.

\begin{Lemma}\label{lem:r}
Let a threshold value $\delta > 0$ and a smooth enough function $f \in L_2(\D,\mu)$ be given. We consider the finite index set $\I_{3\delta}$ such that \eqref{eq:delta_condition} holds. Further, we assume that $3B(\delta_+ + \delta_\Psi) < \min_{\boldk\in\I_{3\delta}} \abs{c_\boldk}$ holds for the detection threshold $\delta_+$ and the projection error threshold $\delta_\Psi$, cf.\ Corollary \ref{cor:single}. Also, let $C \in \R$ be given such that $C \geq C_Q$ for all cubature rules $Q$ used in Algorithm \ref{alg:main}.

Then, choosing the number of detection iterations
\begin{align}\label{eq:detection_iterations}
r \geq \left(1 + \frac{3}{2} B^2 \abs{\I_{3\delta}} + \frac{B^3C}{2\delta} \sum_{\boldh \in \I_{3\delta}^\comp} \abs{c_{\boldh}}\right) (\log 3 + \log d + \log \abs{\I_{3\delta}} - \log \varepsilon)
\end{align}
in Algorithm \ref{alg:main} guarantees that each of the probabilities $\P(E_{1,t}^\comp)$ and $\P(E_{3,t}^\comp)$ is bounded from above by $\varepsilon/(3d)$.
\end{Lemma}
\begin{proof}
We use $\u$ and $K(\I_{3\delta})$ as in Corollary \ref{cor:single} and estimate the probabilities $\P(E_{1,t}^\comp)$ and $\P(E_{3,t}^\comp)$ by \eqref{eq:AtC}. We increase $r$ such that
\begin{align*}
&\P\left( \bigcup_{\boldk \in K(\I_{3\delta})} \left\lbrace \max_{\nu = 1,\hdots,r} \abs{\fhat_{\u,\boldk}^Q (\tilde{\boldx}_\nu)} < \delta_+ \right\rbrace \right) \leq \abs{\I_{3\delta}} \left( \max_{\boldk \in K(\I_{3\delta})} q_{\u,\boldk} \right)^r \leq \frac{\varepsilon}{3d},
\end{align*}
is fulfilled. Consequently, $r$ has to be bounded from below by
\begin{align}\label{eq:r_proof}
\frac{\log 3 + \log d + \log \abs{\I_{3\delta}} - \log \varepsilon}{-\log \max_{\boldk \in K(\I_{3\delta})} q_{\u,\boldk}}.
\end{align}
Hence, we now estimate
\begin{align*}
&\left( \max_{\boldk \in K(\I_{3\delta})} q_{\u,\boldk} \right)^{-1} \\ 
&= \min_{\boldk \in K(\I_{3\delta})} \frac{\displaystyle B \sum_{\boldh=(\boldk,\boldh_{\u^\comp})_\u\in\I_{3\delta}} \abs{c_{\boldh}} + B^2 C_Q \sum_{\boldh \in \I_{3\delta}^\comp} \abs{c_{\boldh}}}{\displaystyle B \sum_{\boldh=(\boldk,\boldh_{\u^\comp})_\u\in\I_{3\delta}} \abs{c_{\boldh}} + B^2 C_Q \sum_{\boldh \in \I_{3\delta}^\comp} \abs{c_{\boldh}} -  B^{-1} \max_{\boldh=(\boldk,\boldh_{\u^\comp})_\u\in\I_{3\delta}} \abs{c_{\boldh}} + \delta_\Psi + \delta_+} \\
&= 1 + \min_{\boldk \in K(\I_{3\delta})} \frac{\displaystyle B^{-1} \max_{\boldh=(\boldk,\boldh_{\u^\comp})_\u\in\I_{3\delta}} \abs{c_{\boldh}} - \delta_\Psi - \delta_+}{\displaystyle B \sum_{\boldh=(\boldk,\boldh_{\u^\comp})_\u\in\I_{3\delta}} \abs{c_{\boldh}} + B^2 C_Q \sum_{\boldh \in \I_{3\delta}^\comp} \abs{c_{\boldh}} -  B^{-1} \max_{\boldh=(\boldk,\boldh_{\u^\comp})_\u\in\I_{3\delta}} \abs{c_{\boldh}} + \delta_\Psi + \delta_+} 
\intertext{and using $3B(\delta_\Psi + \delta_+) < \max_{\boldh=(\boldk,\boldh_{\u^\comp})_\u\in\I_{3\delta}} \abs{c_{\boldh}}$ by assumption yields}
&\geq 1 + \min_{\boldk \in K(\I_{3\delta})} \frac{\displaystyle B^{-1} \max_{\boldh=(\boldk,\boldh_{\u^\comp})_\u\in\I_{3\delta}} \abs{c_{\boldh}} - \delta_\Psi - \delta_+}{\displaystyle B \sum_{\boldh=(\boldk,\boldh_{\u^\comp})_\u\in\I_{3\delta}} \abs{c_{\boldh}} + B^2 C_Q \sum_{\boldh \in \I_{3\delta}^\comp} \abs{c_{\boldh}}} \\
&\geq 1 + \frac{2}{3 B^2} \min_{\boldk \in K(\I_{3\delta})} \frac{\displaystyle \max_{\boldh=(\boldk,\boldh_{\u^\comp})_\u\in\I_{3\delta}} \abs{c_{\boldh}}}{\displaystyle \sum_{\boldh=(\boldk,\boldh_{\u^\comp})_\u\in\I_{3\delta}} \abs{c_{\boldh}} + B C_Q \sum_{\boldh \in \I_{3\delta}^\comp} \abs{c_{\boldh}}}.
\end{align*}
Exploiting the fact that $\max_{\boldh=(\boldk,\boldh_{\u^\comp})_\u\in\I_{3\delta}} \abs{c_{\boldh}} > 3\delta$, the fraction inside the minimum can be bounded from below by
\begin{align*}
& \frac{\displaystyle \max_{\boldh=(\boldk,\boldh_{\u^\comp})_\u\in\I_{3\delta}} \abs{c_{\boldh}}}{\displaystyle \sum_{\boldh=(\boldk,\boldh_{\u^\comp})_\u\in\I_{3\delta}} \abs{c_{\boldh}} + B C_Q \sum_{\boldh \in \I_{3\delta}^\comp} \abs{c_{\boldh}}} \\
\geq & \frac{\displaystyle \max_{\boldh=(\boldk,\boldh_{\u^\comp})_\u\in\I_{3\delta}} \abs{c_{\boldh}}}{\displaystyle \abs{\I_{3\delta}} \max_{\boldh=(\boldk,\boldh_{\u^\comp})_\u\in\I_{3\delta}} \abs{c_{\boldh}} + \frac{\displaystyle \max_{\boldh=(\boldk,\boldh_{\u^\comp})_\u\in\I_{3\delta}} \abs{c_{\boldh}}}{3\delta} B C_Q \sum_{\boldh \in \I_{3\delta}^\comp} \abs{c_{\boldh}}} \\
= & \frac{\displaystyle 1}{\displaystyle \abs{\I_{3\delta}} + \frac{B C_Q}{3 \delta}\sum_{\boldh \in \I_{3\delta}^\comp} \abs{c_{\boldh}}} = \frac{\displaystyle 3 \delta}{\displaystyle 3 \delta \abs{\I_{3\delta}} + B C_Q \sum_{\boldh \in \I_{3\delta}^\comp} \abs{c_{\boldh}}},
\end{align*}
which is now independent of $\boldk$. Hence, we have
\begin{align*}
\left( \max_{\boldk \in K(\I_{3\delta})} q_{\u,\boldk} \right)^{-1} \geq 1 + \frac{\displaystyle 2 \delta B^{-2}}{\displaystyle 3 \delta \abs{\I_{3\delta}} + B C_Q \sum_{\boldh \in \I_{3\delta}^\comp} \abs{c_{\boldh}}}.
\end{align*}
Consequently, and since $\log(x+1)\geq \frac x{x+1}$ for all $x>-1$ and hence $\frac 1{\log(x+1)}\leq \frac {x+1}x = 1+\frac1x$, we obtain for the denominator in \eqref{eq:r_proof} the estimate
\begin{align*}
\frac{1}{\displaystyle -\log \left( \max_{\boldk \in K(\I_{3\delta})} q_{\u,\boldk} \right) } &= \frac{1}{\displaystyle \log \left( \max_{\boldk \in K(\I_{3\delta})} q_{\u,\boldk} \right)^{-1}} \leq \frac{1}{\log \left( 1 + \frac{\displaystyle 2 \delta B^{-2}}{\displaystyle 3 \delta \abs{\I_{3\delta}} + B C_Q \sum_{\boldh \in \I_{3\delta}^\comp} \abs{c_{\boldh}}} \right)} \\
&\leq 1 + \frac{\displaystyle 3 \delta \abs{\I_{3\delta}} + B C_Q \sum_{\boldh \in \I_{3\delta}^\comp} \abs{c_{\boldh}}}{\displaystyle 2 \delta B^{-2}} = 1 + \frac{3}{2} B^2 \abs{\I_{3\delta}} + \frac{B^3C_Q}{2\delta} \sum_{\boldh \in \I_{3\delta}^\comp} \abs{c_{\boldh}}.
\end{align*}
Finally, using $C$ instead of $C_Q$, this bound is applicable independently of $Q$. Therefore, the choice \eqref{eq:detection_iterations} then satisfies the lower bound \eqref{eq:r_proof} and $\P(E_{1,t}^\comp) \leq \varepsilon / (3d)$ and $\P(E_{3,t}^\comp) \leq \varepsilon / (3d)$ are fulfilled.
\end{proof}

\begin{Remark}
The lower bound \eqref{eq:detection_iterations} depends linearly on $\sum_{\boldh \in \I_{3\delta}^\comp} \abs{c_{\boldh}}$, which will be small for reasonable choices of $\delta$ as well as fast enough decaying coefficients $c_\boldh$. Still, we may not access this value directly and need to bound it from above by some more accessible value, e.g.,
\begin{align*}
\sum_{\boldh \in \I_{3\delta}^\comp} \abs{c_{\boldh}} \leq \sum_{\boldh \in \N^d} \abs{c_{\boldh}}
\end{align*}
However, the sum might be tremendously smaller than such bounds, since it does not consider the largest coefficients of $f$, i.e., all basis coefficients $c_\boldh$, whose absolute value is larger than our threshold $3\delta$.

Additionally, note that for sparse functions $f$ the sum $\sum_{\boldh \in \I_{3\delta}^\comp} \abs{c_{\boldh}}$ vanishes completely for a large enough threshold $3\delta$.
\end{Remark}

\subsection{Main result}

Now we collected all ingredients to state and prove our main theorem:

\begin{Theorem}\label{theorem1}
Let a threshold $\delta>0$ and a failure probability $\varepsilon \in (0,1)$ be given. We consider a smooth enough function $f \in L_2(\D,\mu)$ such that the corresponding finite index set $\I_{3\delta}$ fulfilling the condition \eqref{eq:delta_condition} is non-empty. Further, we assume that the projection error terms $\Psi_{\u,\boldk}^{Q,J}(\tilde{\boldx})$ given in \eqref{eq:noise} for $\u = \{t\}, t=1,\ldots,d,$ and $\u = \{1,\ldots,t\}, t=2,\ldots,d-1,$ fulfill the bound
\begin{align}\label{eq:psi_delta}
\norm{\Psi_{\u,\boldk}^{Q,\I_{3\delta}}}_{L_1(\D_{\u^\comp})} < \delta_\Psi
\end{align}
for all $\boldk \in \mathcal{P}_{\u}(\I_{3\delta})$ or $\boldk \in \left( \I^{(1,\hdots,t-1)} \times \I^{(t)} \right) \cap \mathcal{P}_{(1,\hdots,t)}(\I_{3\delta})$, respectively, for some $\delta_\Psi < \delta$ uniformly for all possibly used cubature rules $Q$. Finally, let $C \in \R$ be given such that $C \geq C_Q$ for all possibly used cubature rules $Q$.

We apply Algorithm \ref{alg:main} to the function $f$ using the following parameters. We choose
\begin{itemize}
\item the search space $\Gamma \subset \N^d$ such that $\I_{3\delta} \subset \Gamma$,
\item the detection threshold $\delta_+$ such that $3B(\delta_+ + \delta_\Psi) < \min_{\boldk\in\I_{3\delta}} \abs{c_\boldk}$,
%\item the sparsity parameter $s \geq \abs{\I_{3\delta}}$,
\item the number of detection iterations $ r \geq \left(1 + \frac{3}{2} B^2 \abs{\I_{3\delta}} + \frac{B^3C}{2\delta} \sum_{\boldh \in \I_{3\delta}^\comp} \abs{c_{\boldh}}\right) (\log 3 + \log d + \log \abs{\I_{3\delta}} - \log \varepsilon)$. 
\end{itemize}
We assume that the construction of the cubature rules $Q$ fails with probability at most $\varepsilon/(3d)$.

Then, with a probability $1-\varepsilon$, the index set $\I$ in the output of Algorithm \ref{alg:main} contains the whole index set $\I_{3\delta}$.
\end{Theorem}
\begin{proof}
The assertion follows directly when inserting the probabilities $\P(E_{1,t}^\comp)$ and $\P(E_{3,t}^\comp)$ proven in Lemma \ref{lem:r} in the union bound estimate 
\begin{align*}
\P \left( \bigcap_{t=1}^d E_{1,t} \cap \bigcap_{t=2}^d E_{2,t} \cap \bigcap_{t=2}^d E_{3,t} \right) &= 1 - \P \left( \bigcup_{t=1}^d E_{1,t}^\comp \cup \bigcup_{t=2}^d E_{2,t}^\comp \cup \bigcup_{t=2}^d E_{3,t}^\comp \right) \\
&\geq 1 - \sum_{t=1}^d \P(E_{1,t}^\comp) - \sum_{t=2}^d \P(E_{2,t}^\comp) - \sum_{t=2}^d \P(E_{3,t}^\comp)
\end{align*}
discussed in the beginning of Section \ref{sec:theory}.
\end{proof}

\begin{Remark}
While we already discussed the smoothness assumption in Remark \ref{rem:smooth}, inequality \eqref{eq:psi_delta} adds another technical assumption to Theorem \ref{theorem1}. As mentioned in Section \ref{sec:intro}, $\Psi_{\u,\boldk}^{Q,\I_{3\delta}}$ will vanish for sparse functions $f$ (given a suitable choice of $Q$) and hence fulfill \eqref{eq:psi_delta}. However, also the case of sparsity with additional noise can be covered easily depending on the magnitude and type of the noise. 

For non-sparse functions $f$, we could consider a sufficiently fast decay of the basis coefficients $c_\boldk$. As an example, we refer to the Fourier case and the weigthed Wiener spaces as studied in \cite{PoSchmi19}. Therein, weights $\omega^{\alpha,\beta}(\boldk)$ are considered, which are of product and order-dependent structure and where $\alpha$ regulates the isotropic and $\beta$ the dominating mixed smoothness. In \cite[Sec.~5]{PoSchmi19} the case of scattered data as well as the black-box scenario are investigated, where the latter corresponds to our case here.

Finally, note that a large gap in the size of the basis coefficients will also have a huge effect on the size of $\Psi_{\u,\boldk}^{Q,\I_{3\delta}}$, if $\I_{3\delta}$ is chosen to match that gap.
\end{Remark}

Theorem \ref{theorem1} guarantees the output set $\I$ of Algorithm \ref{alg:main} to contain the index set $\I_{3\delta}$ with probability $1-\varepsilon$. This holds, since we have shown that all the projections $\mathcal{P}_\u(\I_{3\delta})$ are detected with high probability. Note that if the projection error term $\Psi_{\u,\boldk}^{Q,J}$ is large enough in Formula \eqref{eq:proj_coef_theory}, the projected coefficients appear to be large enough and might be detected anyway, even if there is not a single $\boldh \in \I_{3\delta}$ with $\boldh = (\boldk,\boldh_{\u^\comp})_\u$, cf. Remark \ref{rem:falsepos}. This leads to additional index projections detected, which do not belong to the largest coefficients. Theorem \ref{theorem1} implicitely assumes that the cut-off at the end of each step in Algorithm \ref{alg:main} might reduce the amount of such unnecessary detections, but never throws away the important projections of the indices $\boldh \in \I_{3\delta}$. This can be ensured by choosing the sparsity parameters $s$ and $s_{\text{local}}$ in Algorithm \ref{alg:main} large enough. Note that in real applications the optimal choice of those parameters is probably unknown, so they would need to be chosen roughly.

\begin{Remark}\label{rem:falsepos}
We now briefly estimate the probability of additional detections $\boldk \not \in \mathcal{P}_{\u}(\I_{3\delta})$ to get a feeling for how many of those detections we should expect each time. We use the right-hand side inequality
\begin{align*}
\P(\abs{Y} \geq \delta_+) \leq \frac{\displaystyle \E h(Y)}{\displaystyle h(\delta_+)}
\end{align*}
from \cite[Par.\ 9.3.A]{Lo77}, again with $h(t)\coloneqq \abs{t}$ as in the proof of Lemma \ref{lem:trig_prob} and $Y = \Psi_{\u,\boldk}^{Q,J}(\tilde{\boldx})$. Hence, we have the estimate
\begin{align*}
\P\left(\abs{\Psi_{\u,\boldk}^{Q,J}(\tilde{\boldx})} < \delta_+\right) \geq 1-\frac{\norm{\Psi_{\u,\boldk}^{Q,J}}_{L_1(\D_{\u^\comp})}}{\delta_+} > 1-\frac{\delta_\Psi}{\delta_+}
\end{align*}
for the probability that $\Psi_{\u,\boldk}^{Q,J}$ is smaller than the detection threshold $\delta_+$. So we end up with the probability 
\begin{align*}
\P \left(\max_{\nu = 1,\ldots,r} \abs{\Psi_{\u,\boldk}^{Q,J}(\tilde{\boldx}_\nu)} \geq \delta_+\right) < 1- \left( 1-\frac{\delta_\Psi}{\delta_+} \right)^r
\end{align*}
that in at least one detection iteration $\nu=1,\ldots,r$ the projection error term $\Psi_{\u,\boldk}^{Q,J}$ passes the detection threshold $\delta_+$. So for each $\boldk \not \in \mathcal{P}_{\u}(\I_{3\delta})$, there is a chance of at most $1-(1-\frac{\delta_\Psi}{\delta_+})^r$ to be detected accidentally. This matches the intuitive idea that functions with a large $\delta_\Psi$ behave significantly worse than functions with very small coefficients $c_\boldk \ll 3\delta$ for all $\boldk \in \I_{3\delta}^\comp$ or functions, which are nearly sparse.
\end{Remark}

\section{Numerics}
\label{sec:num}

We now investigate the performance and results of our algorithm for two different settings. The first numerical example considers the approximation of a $10$-dimensional, periodic test function in the space $L_2(\T^{10})$ using the Fourier basis. The second part is the approximation of a $9$-dimensional, non-periodic test function in the space $L_2([-1,1]^9,\mu_{\mathrm{Cheb}})$ with $\mu_{\mathrm{Cheb}}$ the Chebyshev measure using the Chebyshev basis, cf. Section \ref{subsec:setup}.

While we mentioned in Section \ref{subsec:method} that an additional recomputation of the basis coefficients on the detected index set in Step 3 of Algorithm \ref{alg:main} is not necessary, the error size of the coefficient approximation might be significantly smaller in this case. To respect this possible lack of accuracy of the herein used version of Algorithm \ref{alg:main}, we investigate the precision of the coefficient approximation and the crucial aim of the algorithm, the detection of a useful sparse index set $\I$, separately.

Finally, note that we do not control the size of the output index set $\I$ by the detection threshold $\delta_+$ as in our approach in Section \ref{sec:theory} in our numerical tests. Here, we follow the more common approach of choosing $\delta_+$ relatively small, but controlling the output using the sparsity parameter $s$. Hence, we do not need to estimate a suitable choice for $\delta_+$ based on the intended threshold $\delta$, cf. Theorem \ref{theorem1}, but also have no theoretical guarantee on the output index set $\I$ anymore.

\subsection{10-dimensional periodic test function}
\label{subsec:ex_four}

For this example, we consider the frequency domain $\Z^d$ instead of $\N^d$, as it is more convenient to use for the Fourier basis. Consider the multivariate test function $f: \T^{10} \rightarrow \R,$
\begin{align*}
f(\boldx) \coloneqq \prod_{j\in\left\lbrace 1,3,8 \right\rbrace} N_2(x_j) + \prod_{j\in\left\lbrace 2,5,6,10 \right\rbrace} N_4(x_j) + \prod_{j\in\left\lbrace 4,7,9 \right\rbrace} N_6(x_j),
\end{align*}
from e.g. \cite[Sec.\ 3.3]{PoVo14} and \cite[Sec.\ 4.2.3]{KaKrVo20}, where $N_m:\T\rightarrow\R$ is the B-Spline of order $m\in\N$,
\begin{align*}
N_m(x) \coloneqq C_m \sum_{k\in\Z} \sinc\left(\frac{\pi}m k \right)^m(-1)^k\e^{2\pi\ii k x},
\end{align*}
with a constant $C_m>0$ such that $\norm{N_m}_{L_2(\T)}=1$. The function $f$ has infinitely many non-zero Fourier coefficients. The largest and therefore most important coefficients are expected to be supported on a union of a three-dimensional symmetric hyperbolic cross in the dimensions 1, 3, 8, a four-dimensional symmetric hyperbolic cross in the dimensions 2, 5, 6, 10, and a three-dimensional symmetric hyperbolic cross in the dimensions 4, 7, 9, each corresponding to the important coefficients of one of the three summands due to their decay properties. Therefore we use the search space $\Gamma \coloneqq \tilde{H}_{N}^{10,\frac12}$, where 
\begin{align}\label{eq:sym_hc}
\tilde{H}_{N}^{d,\gamma} \coloneqq \left\lbrace \boldk \in \Z^d: \prod_{j=1}^d \max \left(1,\frac{\abs{k_j}}{\gamma}\right) \leq 2^N \right\rbrace 
\end{align}
is the $d$-dimensional symmetric hyperbolic cross with weight $0<\gamma\leq 1$ and extension $N>0$ such that $2^N\in\N$. Moreover, we set the number of detection iterations $r = 5$. We increase the sparsity $s$ exponentially and use the local sparsity parameter $s_\mathrm{local} = 1.2 s$. Due to the fast decay of the Fourier coefficients and the increasing sparsity $s$, we fix the detection threshold $\delta_+ = 10^{-12}$.

The algorithm was implemented in MATLAB\textsuperscript{\textregistered} and tested using 2 six core CPUs Intel\textsuperscript{\textregistered} Xeon\textsuperscript{\textregistered} CPU E5-2620 v3 @ 2.40GHz and 64 GB RAM. All tests are performed 10 times and the relative $L_2(\T^{10})$ approximation error
\begin{align*}
\frac{\norm{f-S_{\I}f}_{L_2(\T^{10})}}{\norm{f}_{L_2(\T^{10})}} \coloneqq \frac{\sqrt{\norm{f}_{L_2(\T^{10})}^2 - \sum_{\boldk\in\I} \vert c_{\boldk}\vert^2}}{\norm{f}_{L_2(\T^{10})}}
\end{align*}
as well as the coefficient approximation error
\begin{align*}
\norm{\boldsymbol{\fhat}-\boldsymbol{c}}_{\ell_p(\I)} \coloneqq \begin{cases}
\left( \sum_{\boldk\in\I} \vert\fhat_\boldk - c_{\boldk}\vert^p \right)^{\frac{1}{p}} & 1\leq p < \infty \\
\max_{\boldk\in\I} \vert\fhat_\boldk - c_{\boldk}\vert &  p = \infty
\end{cases}
\end{align*}
for $p=1,2,\infty$ are computed. Note that the relative $L_2(\T^{10})$ approximation error uses the exact coefficients $c_{\boldk}$ and hence only depends on the computed index set $\I$. Therefore, this error indicates how well our detected indices $\boldk \in \I$ can approximate the function $f$ without possible loss due to non-optimal approximations $\fhat_\boldk$ of the coefficients $c_{\boldk}$.

We consider three different cubature methods $Q$ to construct and evaluate the approximated projected (Fourier) coefficients $\fhat_{\u,\boldk}^Q (\tilde{\boldx})$ in Algorithm \ref{alg:main}:

\begin{itemize}
\item \textit{Monte Carlo points (MC)}: We draw $M \coloneqq \abs{K}\log(\abs{K})$ nodes $\boldxi_j, j=1,\ldots,M,$ uniformly at random in $\T^t$ and set $w_j=M^{-1}$ for all $j=1,\ldots,M$. To improve the accuracy, we subsequent apply the least squares method with up to $20$ iterations and the tolerance $10^{-6}$. Hence, the method is no longer an equally weighted Monte Carlo cubature.

\item \textit{single rank-1 lattices (R1L)}:
We construct a rank-1 lattice which is a spatial discetization for $K$ using \cite[Algorithm~5]{Kae23}. An FFT approach is used in order to compute the projected coefficients simultaneously, cf.\ \cite{LiHi03} and \cite[Algorithm~3.2]{kaemmererdiss}, which is equivalent to the application of a cubature rule using the $M$ sampling nodes of the rank-1 lattice and the weights $w_j=M^{-1}$.

\item \textit{multiple rank-1 lattices (MR1L)}:
We construct a spatial discretization for $K$ using \cite[Algorithm~5]{Kae17}. The discretization consists of a set of rank-1 lattices whose structure allows efficient computations of the Fourier matrix and its adjoint matrix, cf.\ \cite{Kae16}, which in turn is used to apply the least squares method to compute the projected Fourier coefficients. The calculation is equivalent to applying a non-equally weighted cubature rule to calculate each projected coefficient.
\end{itemize}

Table \ref{tbl:comp1} states upper bounds on the sampling and computational complexity of the methods using the rough bounds $\OO{drS_Q(r^2s^2,d)}$ and $\OO{drT_Q(r^2s^2,d)}$ from Section \ref{subsec:comp}. Those complexities are simplified using particular assumptions. For more detail on the complexities of the reconstruction methods, see the respective references.

\renewcommand{\arraystretch}{1.1}
\begin{table}[tb]
\centering
\begin{tabular}{c|c|c|c} 
 & sampling complexity & computational complexity & comments \\ \hline
 MC & $\OO{d r^3 s^2 \log(rs)}$ & $\OO{d r^5 s^4 \log(rs)}$ & - \\ \hline
 R1L & $\OO{d r^5 s^4}$ & $\OO{(r^2 s^2 + d \log^{1+\varepsilon}(d))d r^3 s^3 \log(rs)}$ & (a) \\ \hline
 MR1L & $\OO{d r^3 s^2 \log(rs)}$ & $\OO{d^2r^3s^2 \log^3(rs)}$ & (b) \\
\end{tabular}
\caption{Rough bounds on the sampling and computational complexity of Algorithm \ref{alg:main} for the three different reconstruction methods using the results from Section \ref{subsec:comp}. Comments: (a) We are assuming $r^2s^2 > \max_{\boldk \in \Gamma}\vert\vert \boldk \vert\vert_\infty$. The computational complexity holds with high probability $1-\delta$ assuming $\varepsilon>0$ and fixed failure probability $\delta$. (b) This holds with high probability assuming fixed failure probability $\delta$ and oversampling factor $c$. The computational complexity is assuming $r^2 s^2 \log^2(rs) > (\max_{\boldk \in \Gamma}\vert\vert \boldk \vert\vert_\infty)^2$ and a fixed upper bound on the amount of least squares iterations.}
\label{tbl:comp1}
\end{table}
\renewcommand{\arraystretch}{1}

\begin{figure}
\begin{subfigure}[c]{0.5\textwidth}
	\begin{center}
	\includegraphics{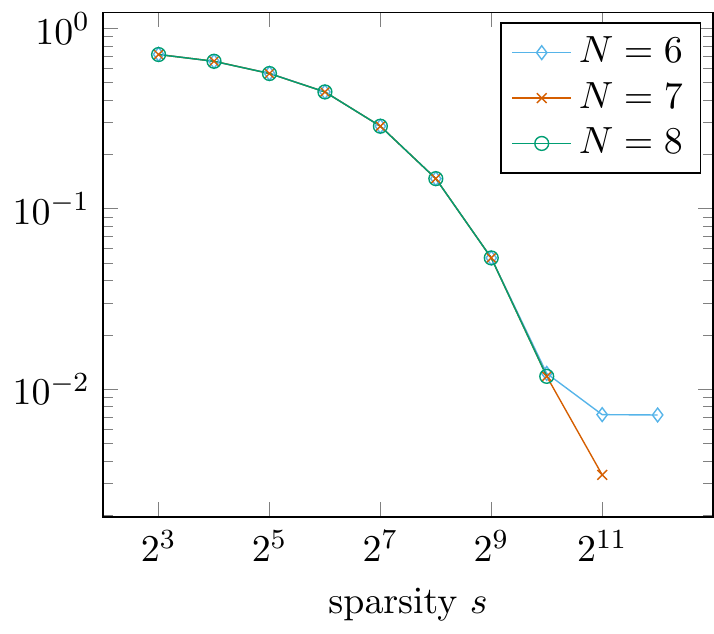}	
	\end{center}
	\vspace{-1\baselineskip}
	\caption{relative $L_2(\T^{10})$ approx. error for MC}\label{fig:MC10d}
	\vspace{+1\baselineskip}
\end{subfigure}
~
\begin{subfigure}[c]{0.5\textwidth}
	\begin{center}
	\includegraphics{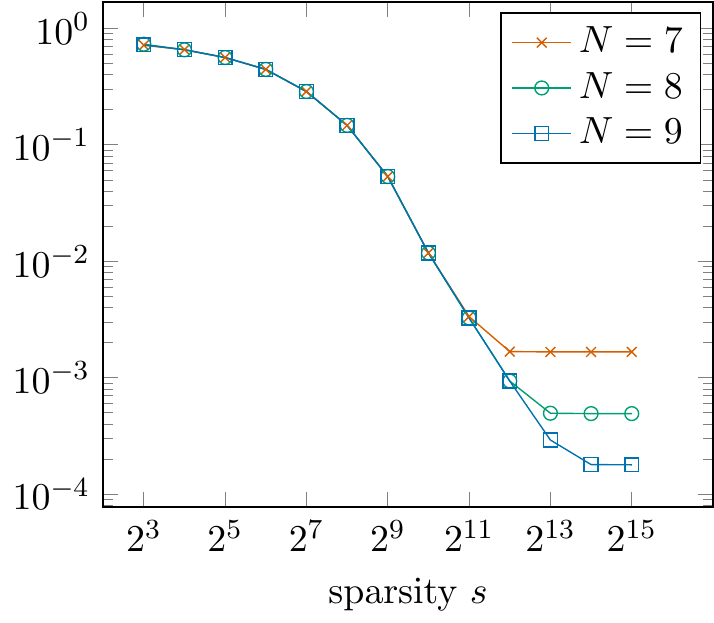}	
	\end{center}
	\vspace{-1\baselineskip}
	\caption{relative $L_2(\T^{10})$ approx. error for R1L}
	\vspace{+1\baselineskip}
\end{subfigure}
~
\vspace*{-0.05cm}
\begin{subfigure}[c]{0.5\textwidth}
	\begin{center}
	\includegraphics{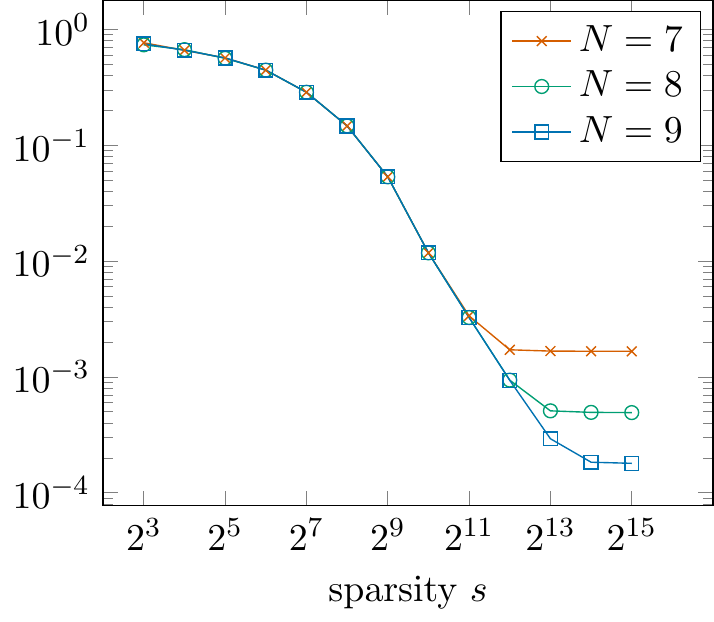}	
	\end{center}
	\vspace{-1\baselineskip}
	\caption{relative $L_2(\T^{10})$ approx. error for MR1L}
	\vspace{+1\baselineskip}
\end{subfigure}
~
\begin{subfigure}[c]{0.5\textwidth}
	\begin{center}
	\includegraphics{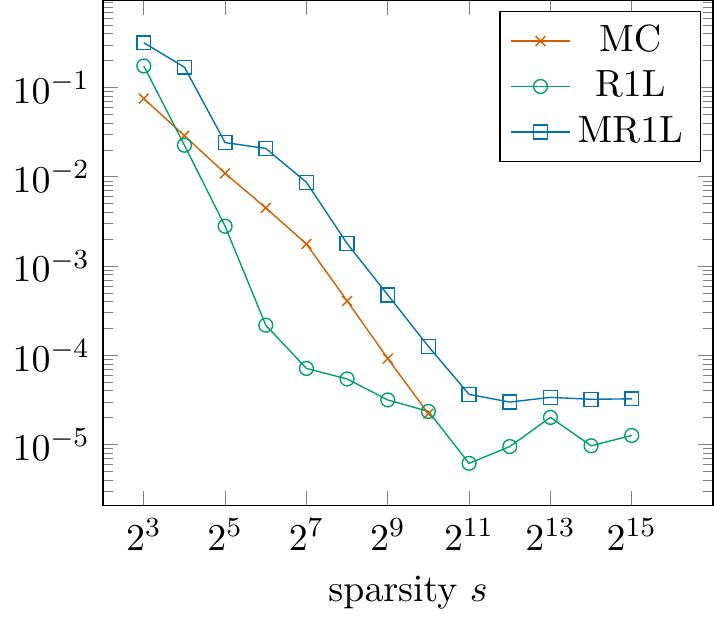}	
	\end{center}
	\vspace{-1\baselineskip}
	\caption{coef. approx. error in $\ell_\infty$ for $N=8$}
	\vspace{+1\baselineskip}
	\label{fig:err10d}
\end{subfigure}
~
\vspace*{-0.05cm}
\begin{subfigure}[c]{0.5\textwidth}
	\begin{center}
	\includegraphics{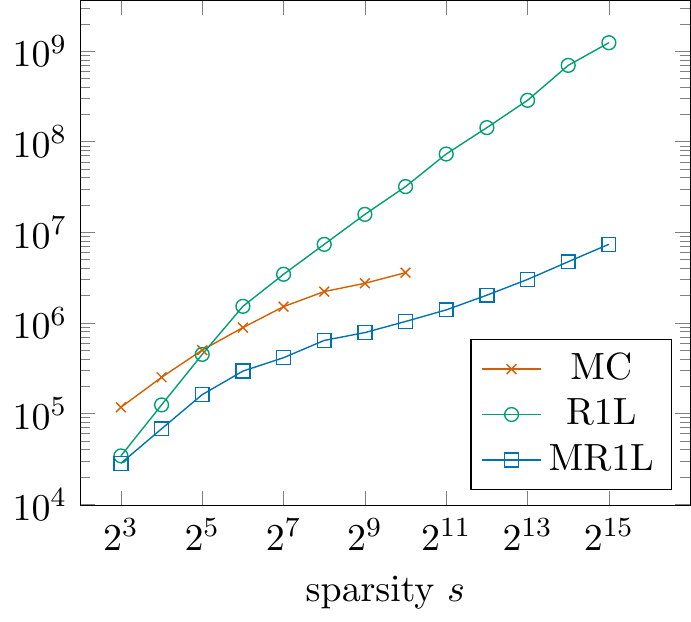}	
	\end{center}
	\vspace{-1\baselineskip}
	\caption{amount of samples for $N=8$}\label{fig:sam10d}
	\vspace{+1\baselineskip}
\end{subfigure}
~
\begin{subfigure}[c]{0.5\textwidth}
	\begin{center}
	\includegraphics{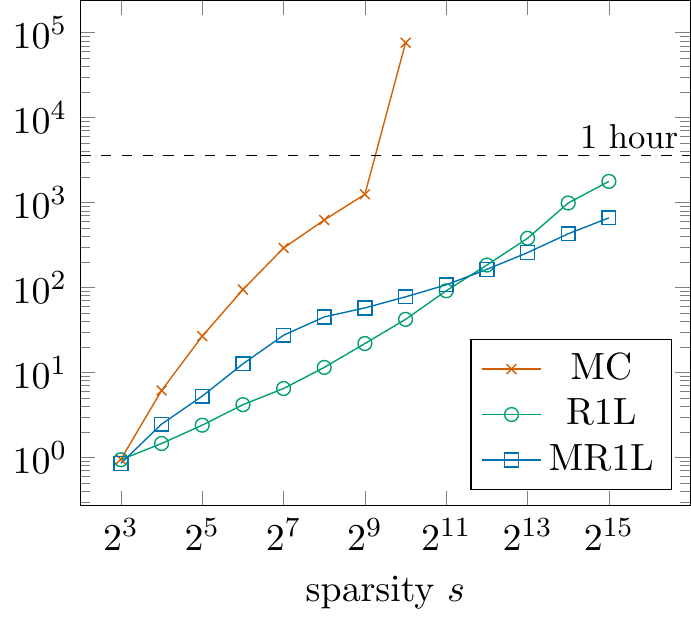}	
	\end{center}
	\vspace{-1\baselineskip}
	\caption{computation time (in seconds) for $N=8$}
	\vspace{+1\baselineskip}
\end{subfigure}
\vspace*{-0.7cm}
\caption{Approximation results for the $10$-dimensional periodic test function}
\label{fig:test10d}
\end{figure}

Figure \ref{fig:test10d} illustrates the decay of the relative $L_2(\T^{10})$ approximation error for the three different cubature methods used as well as the $\ell_\infty$ coefficient approximation error, the amount of samples and the computation time. Note that these are the medians of the 10 test runs. While we also performed tests for the rank-1 lattice approaches for higher extensions up to $N=13$ and larger sparsities like $s=2^{17}$, the MC approach failed to deliver already for smaller parameters because of the large matrices used there. Figure \ref{fig:MC10d} only shows the data of those tests that could be performed successfully. The computation time also increases at a much faster rate, while all illustrated tests for R1L and MR1L computed in less than an hour. Anyways, the detected index sets $\I$ seem to perform very well for all three methods and larger extensions $N$ allow an even longer decay of the relative $L_2(\T^{10})$ approximation error w.r.t.\ the sparsity $s$. 

The $\ell_\infty$ coefficient approximation error decays well for all examples, but starts to stagnate a bit earlier w.r.t.\ the sparsity than the relative error. This effect seems to be caused by the aliasing of the coefficients surrounding $\Gamma$, since for larger or smaller values of $N$ the stagnation also starts later or earlier, respectively. Analogously, the $\ell_1$ and $\ell_2$ coefficient approximation error decay in a similar way but are not illustrated here to preserve clarity. The amount of samples illustrated in Figure \ref{fig:sam10d} grows reasonably. As expected, the R1L approach needs most samples while the growth for the MR1L and MC approaches is considerably slower. The amount of samples tends to be very large compared to the size of our search space $\Gamma$, e.g., $\vert\tilde{H}_{8}^{10,\frac12}\vert \approx 2.39 \cdot 10^6$, but other hyperbolic crosses or even the full grid should result in comparable amounts of samples as long as the extension $N$ is of similar size, while the amount of candidates grows tremendously for those search spaces.

\subsection{9-dimensional non-periodic test function}
\label{subsec:ex_cheb}

\begin{figure}
\begin{subfigure}[c]{0.5\textwidth}
	\begin{center}
		\includegraphics{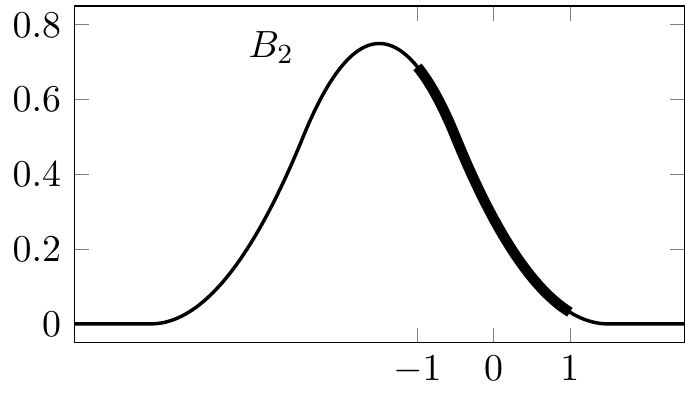}		
	\end{center}
\end{subfigure}
~
\begin{subfigure}[c]{0.5\textwidth}
	\begin{center}
\includegraphics{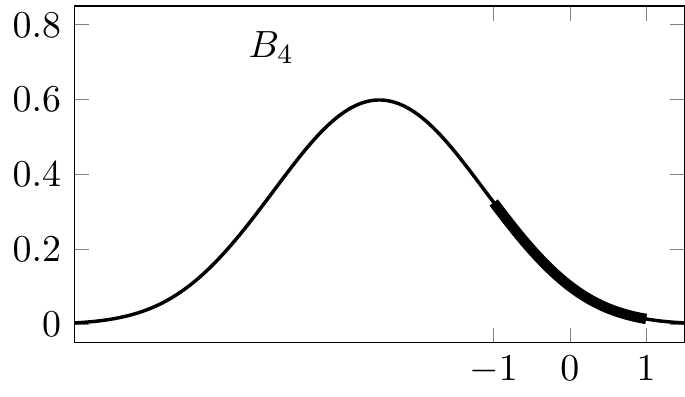}	
	\end{center}
\end{subfigure}
\caption{B-Splines $B_2$ and $B_4$ and the considered domain $[-1,1]$}
\label{fig:bspline}
\end{figure}

Consider the multivariate test function $f: [-1,1]^{9} \rightarrow \R,$
\begin{align*}
f(\boldx) \coloneqq \prod_{j\in\left\lbrace 1,3,4,7 \right\rbrace} B_2(x_j) + \prod_{j\in\left\lbrace 2,5,6,8,9 \right\rbrace} B_4(x_j),
\end{align*}
from e.g. \cite[Sec.\ 4.2.2]{volkmerdiss} and \cite[Sec.\ III.B]{PoVo17}, where $B_2:\R\rightarrow\R$ and $B_4:\R\rightarrow\R$ are shifted, scaled and dilated B-Splines of order $2$ and $4$, respectively, see Figure \ref{fig:bspline} for an illustration and \cite{volkmerdiss} for their rigorous definition.

As in Section \ref{subsec:ex_four}, the function $f$ is not sparse in the Chebyshev frequency domain and we expect the significant Chebyshev coefficients $c_{\boldk}, \boldk \in \N^9,$ of the function to be supported on a union of a four-dimensional hyperbolic cross like structure in the dimensions $1, 3, 4, 7,$ and a five-dimensional hyperbolic cross like structure in the dimensions $2, 5, 6, 8, 9$. Hence, we restrict ourselves to the search space $\Gamma \coloneqq \tilde{H}_n^{9,\frac12} \cap \N^9$ with $\tilde{H}_{n}^{d,\gamma}$ as defined in \eqref{eq:sym_hc}. Again, we increase the sparsity $s$ while fixing the number of detection iterations $r=5$, the detection threshold $\delta_+=10^{-12}$ and the local sparsity parameter $s_{\mathrm{local}}=1.2 s$.

All tests are performed 10 times and the relative $L_2([-1,1]^9,\mu_{\mathrm{Cheb}})$ approximation error
\begin{align*}
\frac{\norm{f-S_{\I}f}_{L_2([-1,1]^9,\mu_{\mathrm{Cheb}})}}{\norm{f}_{L_2([-1,1]^9,\mu_{\mathrm{Cheb}})}} \coloneqq \frac{\sqrt{\norm{f}_{L_2([-1,1]^9,\mu_{\mathrm{Cheb}})}^2 - \sum_{\boldk\in\I} \vert c_{\boldk}\vert^2}}{\norm{f}_{L_2([-1,1]^9,\mu_{\mathrm{Cheb}})}}
\end{align*}
as well as the coefficient approximation error
\begin{align*}
\norm{\boldsymbol{\fhat}-\boldsymbol{c}}_{\ell_p(\I)} \coloneqq \begin{cases}
\left( \sum_{\boldk\in\I} \vert\fhat_\boldk - c_{\boldk}\vert^p \right)^{\frac{1}{p}} & 1\leq p < \infty \\
\max_{\boldk\in\I} \vert\fhat_\boldk - c_{\boldk}\vert &  p = \infty
\end{cases}
\end{align*}
for $p=1,2,\infty$ are computed. As before, the relative $L_2([-1,1]^9,\mu_{\mathrm{Cheb}})$ approximation error only depends on the detected index set $\I$ and not on the approximated coefficients $\fhat_\boldk, \boldk \in \I,$ cf.\ Section \ref{subsec:ex_four}.

We consider two different cubature methods $Q$ to construct and evaluate the approximated projected (Chebyshev) coefficients $\fhat_{\u,\boldk}^Q (\tilde{\boldx})$ in Algorithm \ref{alg:main}: \begin{itemize}
\item \textit{Monte Carlo points (cMC)}: We set $m=\sum_{\boldk\in K} 2^{\norm{\boldk}_0}$, cf. \cite[Sec.~1.2]{KuoMiNoNu19}, and draw $M \coloneqq m \log(m)$ nodes $\boldxi_j, j=1,\ldots,M,$ in $[-1,1]$ at random w.r.t.\ the Chebyshev measure $\mu_{\mathrm{Cheb}}$ and set $w_j=M^{-1}$ for all $j=1,\ldots,M$. To improve the accuracy, we again apply the least squares method with up to $20$ iterations and the tolerance $10^{-6}$. Again, the method is now no longer an equally weighted Monte Carlo cubature in the classical sense.
\item \textit{Chebyshev multiple rank-1 lattices (cMR1L)}: Similar to the multiple rank-1 lattice approach in Section~\ref{subsec:ex_four}, there exists a strategy \cite{Kae23_cheb} for discretizing spans of multivariate Chebyshev polynomials using sets of transformed rank-1 lattices \cite{KuoMiNoNu19}. The computation of the evaluation matrix as well as its adjoint can be implemented in an efficient way. We compute the Chebyshev coefficients using the least squares method which is thus equivalent to applying a non-equally weighted cubature rule to calculate each projected coefficient.
\end{itemize}

Table \ref{tbl:comp2} again states upper bounds on the sampling and computational complexity of those methods using the rough bounds $\OO{drS_Q(r^2s^2,d)}$ and $\OO{drT_Q(r^2s^2,d)}$ from Section \ref{subsec:comp}. As before, those complexities are simplified using particular assumptions, while the detailed versions can be found in the given references.

\renewcommand{\arraystretch}{1.1}
\begin{table}[tb]
\centering
\begin{tabular}{c|c|c|c} 
 & sampling complexity & computational complexity & comments \\ \hline
 cMC & $\OO{d r m_\ast \log(m_\ast)}$ & $\OO{d r m_\ast^2 \log(m_\ast)}$ & (a) \\ \hline
 cMR1L & $\OO{d r m_\ast \log(m_\ast)}$ & $\OO{d^2 rm_\ast \log^3(m_\ast)}$ & (a), (b) \\
\end{tabular}
\caption{Rough bounds on the sampling and computational complexity of Algorithm \ref{alg:main} for the two different reconstruction methods using the results from Section \ref{subsec:comp}. Comments: (a) In the worst case, there holds the very pessimistic bound $m_\ast \leq r^2 s^2 2^d$. However, if $\tilde{d}<d$, cf. Section \ref{subsec:gamma}, there holds $m_\ast \leq r^2 s^2 2^{\tilde{d}}$. Other assumptions may yield even better and more realistic bounds. As an example, see \cite[Lem.~1]{KuoMiNoNu19}, where downward closedness yields $m_\ast \leq (r^2 s^2)^{\ln 3/\ln 2}$. (b) This holds with high probability assuming fixed failure probability $\delta$ and oversampling factor $c$. The computational complexity is assuming $r^2 s^2 \log^2(rs) > (\max_{\boldk \in \Gamma}\vert\vert \boldk \vert\vert_\infty)^2$ and a fixed upper bound on the amount of least squares iterations.}
\label{tbl:comp2}
\end{table}
\renewcommand{\arraystretch}{1}

\begin{figure}
\begin{subfigure}[c]{0.5\textwidth}
	\begin{center}
\includegraphics{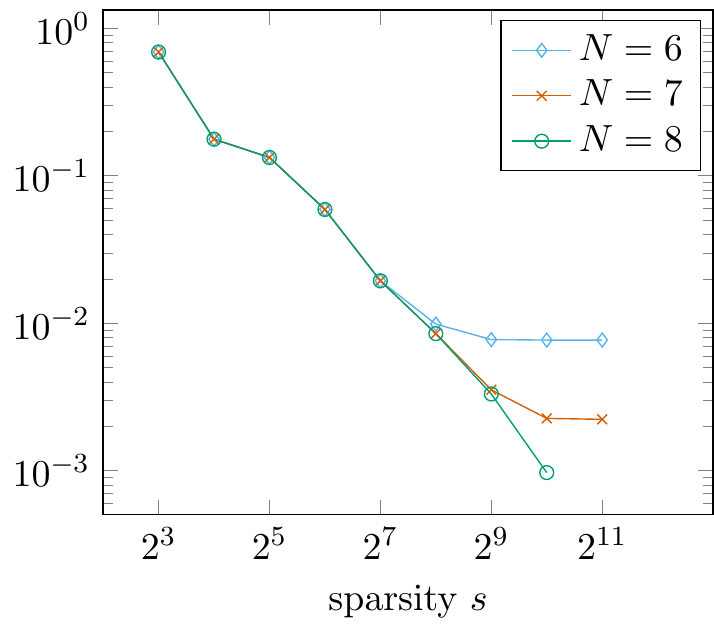}
	\end{center}
	\vspace{-1\baselineskip}
	\caption{relative $L_2([-1,1]^9,\mu_{\mathrm{Cheb}})$ approx. error for cMC}
	\vspace{+1\baselineskip}
\end{subfigure}
~
\begin{subfigure}[c]{0.5\textwidth}
	\begin{center}
\includegraphics{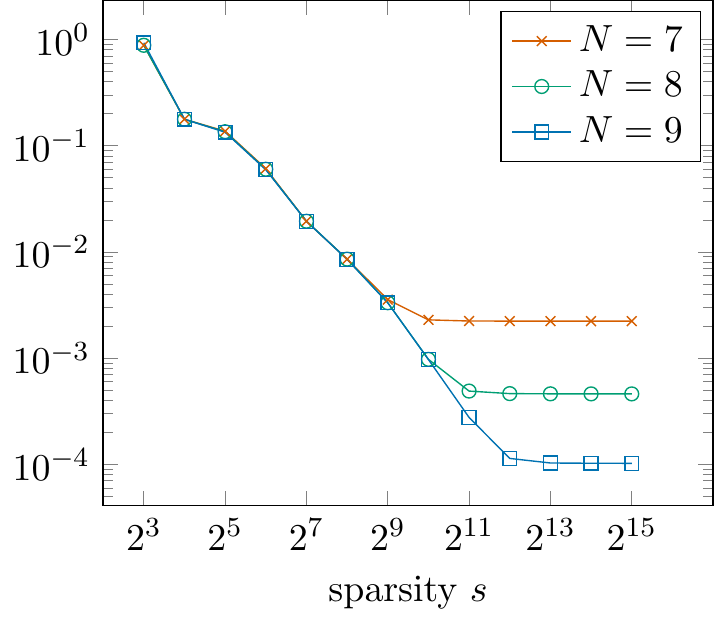}
	\end{center}
	\vspace{-1\baselineskip}
	\caption{relative $L_2([-1,1]^9,\mu_{\mathrm{Cheb}})$ approx. error for cMR1L}
	\vspace{+1\baselineskip}
\end{subfigure}
~
\begin{subfigure}[c]{0.5\textwidth}
	\begin{center}
\includegraphics{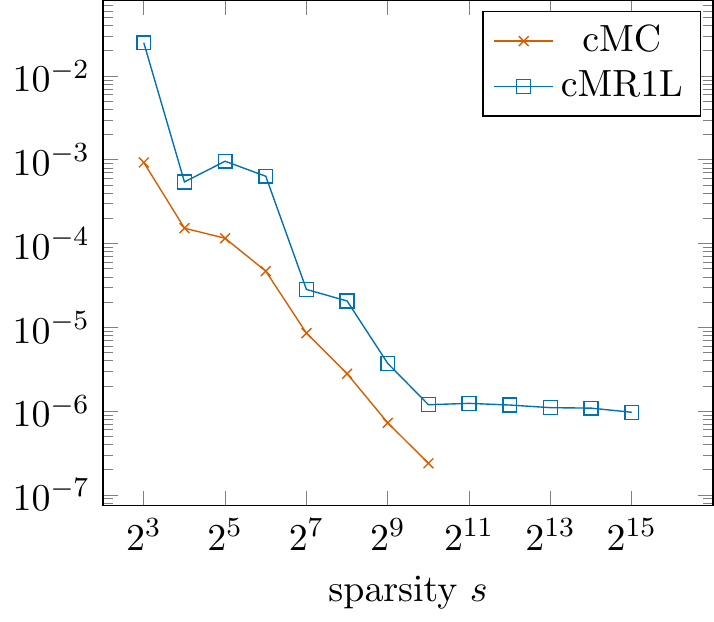}
	\end{center}
	\vspace{-1\baselineskip}
	\caption{coef. approx. error in $\ell_\infty$ for $N=8$}
	\vspace{+1\baselineskip}
\end{subfigure}
~
\begin{subfigure}[c]{0.5\textwidth}
	\begin{center}
\includegraphics{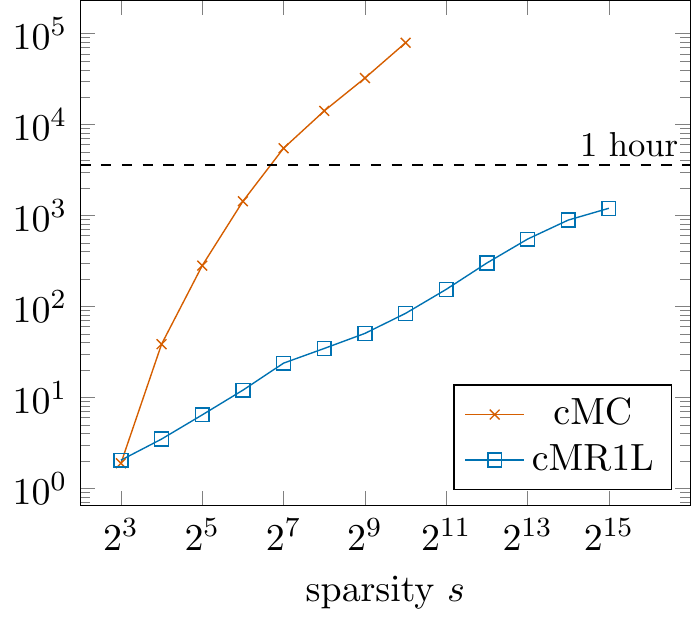}
	\end{center}
	\vspace{-1\baselineskip}
	\caption{computation time (in seconds) for $N=8$}
	\vspace{+1\baselineskip}
\end{subfigure}
\vspace*{-0.7cm}
\caption{Approximation results for the $9$-dimensional non-periodic test function}
\label{fig:test9d}
\end{figure}

Figure \ref{fig:test9d} again illustrates the decay of the relative $L_2([-1,1]^9,\mu_{\mathrm{Cheb}})$ approximation error for the two different methods used as well as the $\ell_\infty$ coefficient approximation error and the computation time. Again, the $\ell_1$ and $\ell_2$ coefficient approximation error decay in a similar way as the $\ell_\infty$ error and are not shown here. As in Section \ref{subsec:ex_four}, we could not apply the cMC method for larger parameters $s$ and $N$ due to the higher computation time as well as the larger amount of samples needed. On the other hand, the cMR1L method transferred their efficiency to the whole dimension-incremental method, managing a significantly shorter computation time and less samples needed. While both the relative $L_2([-1,1]^9,\mu_{\mathrm{Cheb}})$ approximation error as well as the $\ell_\infty$ coefficient approximation error decay as expected for both approaches, the coefficient approximation error is again a little larger for the lattice approach. Those results again underline the importance of the efficiency and accuracy of the underlying reconstruction method to our algorithm. 

\subsection{Reliability}

\begin{figure}
\begin{subfigure}[c]{0.5\textwidth}
\includegraphics{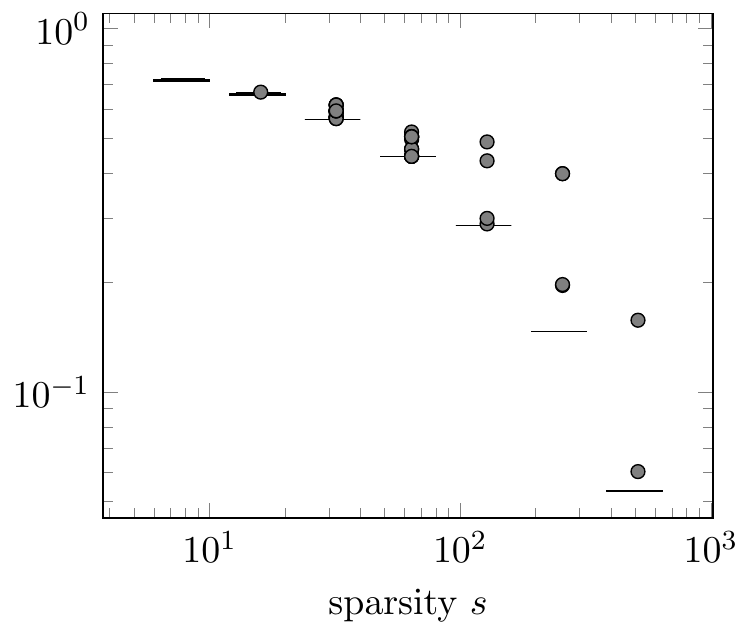}
	\vspace{-0.5\baselineskip}
	\caption{relative approx. error for MC with $N=7$}
\end{subfigure}
~
\begin{subfigure}[c]{0.5\textwidth}
\includegraphics{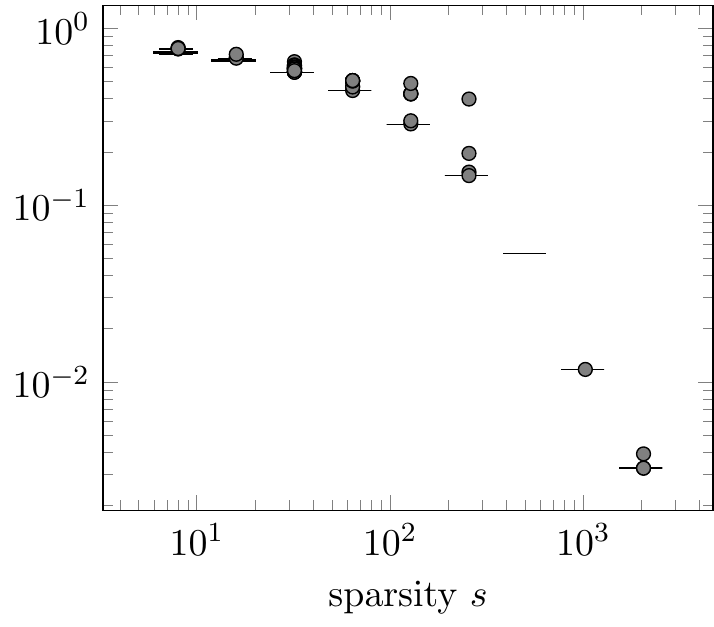}
	\vspace{-0.5\baselineskip}
	\caption{relative approx. error for R1L with $N=8$}
\end{subfigure}
~
\begin{subfigure}[c]{0.5\textwidth}
\includegraphics{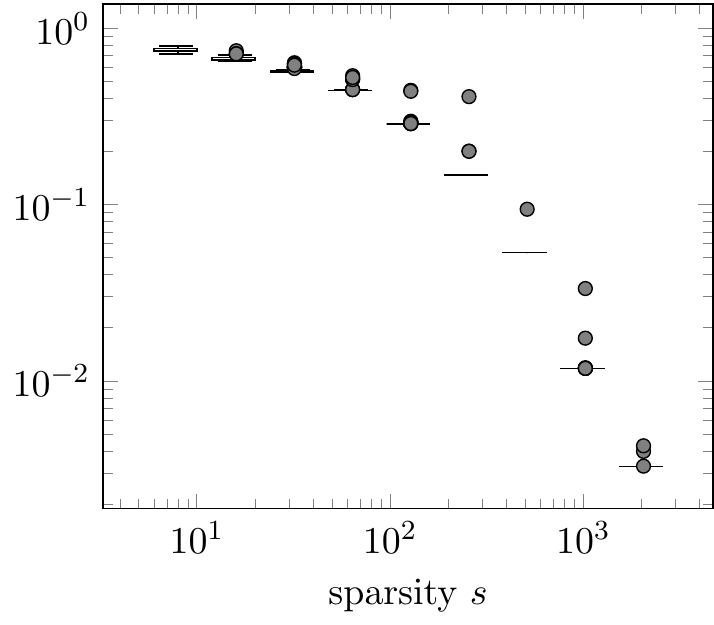}
	\vspace{-0.5\baselineskip}
	\caption{relative approx. error for MR1L with $N=8$}
\end{subfigure}
~
\begin{subfigure}[c]{0.5\textwidth}
\includegraphics{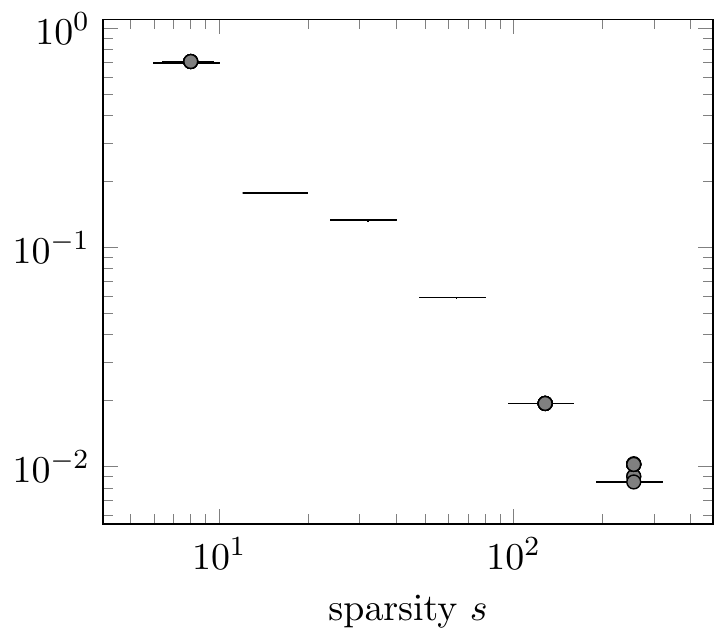}
	\vspace{-0.5\baselineskip}
	\caption{relative approx. error for cMC with $N=7$}
\end{subfigure}
~
\begin{subfigure}[c]{0.5\textwidth}
\includegraphics{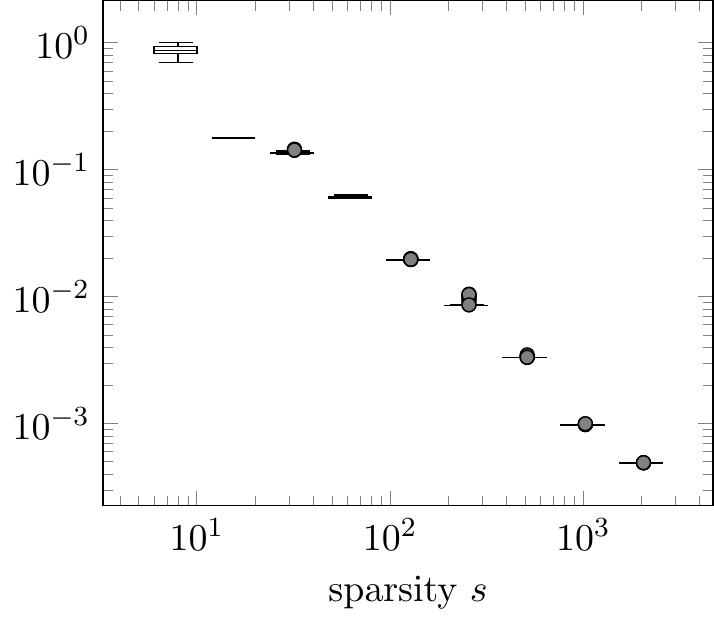}
	\vspace{-0.5\baselineskip}
	\caption{relative approx. error for cMR1L with $N=8$}
\end{subfigure}
\caption{Box-and-whisker plots with outliers for the relative $L_2(\T^{10},\mu)$ and $L_2([-1,1]^9,\mu_{\mathrm{Cheb}})$ approximation errors, respectively, for the different methods with $100$ performed tests per sparsity each.}
\label{fig:error_pics}
\end{figure}

The numerical experiments in Sections \ref{subsec:ex_four} and \ref{subsec:ex_cheb} are regulated by the sparsity $s$ instead of the detection threshold $\delta_+$. This approach seems more natural for applications of our algorithm, but is obviously missing a proper theoretical foundation like Theorem \ref{theorem1}. However, the results already seem very promising in the previous sections. To underline this, we further investigate the reliability of Algorithm \ref{alg:main} by performing chosen tests $100$ instead of $10$ times and plotting the results for the respective relative $L_2$ approximation error in box-and-whisker plots in Figure \ref{fig:error_pics}. We stick to the typical approach by classifying results as outliers, if they are more than 1.5 interquartile ranges above the upper quartile.

The main observation is how small the boxes (and the extension of their whiskers) are for almost all tests, which indicates a high reliability of our algorithm. The amount of outliers is below $10\%$ in most tests, the highest detect amount is $24\%$. However, most of those outliers are still incredibly close to the average results and are most likely to be outliers only due to the tremendously small interquartile ranges. As can be seen in Figure \ref{fig:error_pics}, there are only singular outliers where the accuracy is considerably worse than expected. 

Surprisingly, our approaches in the Chebyshev setting (cMC and cMR1L) seem to behave even better in terms of reliability, as there are even less and also less bad outliers. However, it is not entirely clear whether this is caused by the particular basis or due to side-effects.

\section{Conclusion}
\label{sec:summary}

The presented algorithm is capable of approximating high-dimensional functions very well by detecting a sparse truncation of its basis expansion in the corresponding space. Given a suitable coefficient reconstruction method like rank-1 lattice approaches, the highly adaptive algorithm can be applied to any bounded orthonormal product basis and also benefits tremendously, if the reconstruction method is efficient in terms of, e.g., sample complexity, memory complexity or computational complexity. If several reconstruction methods are available, one should prioritize those with the best properties for the considered situation, e.g., a sample efficient method if sampling is very expensive compared to the rest of the algorithm as in \cite{KaPoTa22}. 

We provide a theoretical reconstruction guarantee for a special kind of methods, which can be seen as a blueprint for similar proofs for other reconstruction methods. On the other hand, Theorem \ref{theorem1} also brings up various open questions which are suitable for further research. It is still unknown how to properly include the sparsity $s$ as cut-off parameter, which is way more suitable for regulating applications of the algorithm, into the theoretical results instead of the detection threshold $\delta_+$. Also, improved bounds on the number of detection iterations $r$ are still desired since the working choice of a small and constant amount does not coincide with the theoretical bound.

Our numerical tests result in promising and reliable nonlinear approximations for the well-known Fourier case as well as the non-periodic Chebyshev case. These results strongly motivate applying the dimension-incremental algorithm to several other bounded orthonormal product bases as in \cite{ChIwVo21}.

Finally, we stated several modifications and improvements of the algorithm throughout the paper, which should be considered in future works to increase the power of the dimension-incremental method even further.

\section*{Acknowledgement}

L. K\"ammerer gratefully acknowledges funding by the Deutsche \linebreak Forschungsgemeinschaft (DFG, German Research Foundation) with the project number \linebreak 380648269 and Daniel Potts with the project number 416228727 -- SFB 1410.

%%%%%%%%%%%%%%%%%%%%%%%%%%%%%%%%%%%%%%%%%%%%%%%%%%%%%%%%%%%%%%%%%%%%%%%%%%%%%%

\bibliographystyle{abbrv}

\end{document}